\documentclass[a4paper]{amsart}

\usepackage{color}
\makeindex
\usepackage{graphicx}        
\usepackage[all]{xy}
\usepackage{amssymb}
\usepackage{amsmath}
\usepackage{amsthm}
\usepackage{hhline}
\usepackage{tikz-cd}
\usepackage{longtable}
\usepackage{ulem}
\usepackage[a4paper, left=2cm, right=2cm, top=2cm,bottom=2cm]{geometry}
\usepackage{hyperref}
\usepackage[draft=true]{minted}
\usepackage{xspace}
\usepackage[backend=bibtex, bibencoding=utf8, doi=false, url=false, isbn=false]{biblatex}
\theoremstyle{definition}
\newtheorem{theorem}{Theorem}[section]
\newtheorem{definition}[theorem]{Definition}
\newtheorem{remark}[theorem]{Remark}

\newtheorem{proposition}[theorem]{Proposition}
\newtheorem{lemma}[theorem]{Lemma}
\newtheorem{corollary}[theorem]{Corollary}

\usepackage[all]{xy}
\usepackage{pdflscape}
\newcommand{\RR}{\mathbb{R}}

\newcommand{\PP}{\mathbb{P}}
\newcommand{\IA}{\mathbb{A}}
\newcommand{\FF}{\mathbb{F}}

\newcommand{\NN}{\mathbb{N}}
\newcommand{\Spec}{\mathrm{Spec}\,}
\newcommand{\CC}{\mathbb C}

\newcommand{\QQ}{\mathbb Q}
\newcommand{\bbk}{\mathbf k}
\newcommand{\ZZ}{\mathbb Z}
\renewcommand{\PP}{\mathbb P}
\newcommand{\GL}{\mathrm{GL}}

\newcommand{\OO}{\mathcal{O}}

\newcommand{\Aut}{\mathrm{Aut}}
\newcommand{\Triv}{\mathrm{Triv}}

\newcommand{\rank}{\mathrm{rank}\,}

\newcommand{\NS}{\mathrm{NS}}
\newcommand{\Num}{\mathrm{Num}}
\newcommand{\Pic}{\mathrm{Pic}}
\newcommand{\Oscar}{OSCAR\xspace}
\newcommand{\D}{\operatorname{d}\!}

\begin{document}

\title{An explicit Enriques surface with an automorphism of minimum entropy}

\author{Simon Brandhorst}
\address{
  Simon Brandhorst: 
  Universit\"at des Saarlandes,
  66123 Saarbr\"ucken,
  Germany 
}
%
\author{Matthias Zach}
\address{
  Matthias Zach:
  RPTU Kaiserslautern-Landau,
  67663 Kaiserslautern,
  Germany 
}

\begin{abstract}
  We derive explicit equations for the Oguiso-Yu automorphism of minimum topological entropy on a complex Enriques surface.
The approach is computer aided and makes use of elliptic fibrations.
\end{abstract}
\maketitle

\section{Introduction}

Let $X$ be a smooth, connected, compact complex K\"ahler surface and $f \in \Aut(X)$ be an automorphism. The topological entropy $h(f)$ of $f$ measures the dynamical complexity of $f$. By theorems of Gromov and Yomdin \cite{Yomdin87,gromov90} it is the logarithm
\[h(f)=\log d_1(f)\]
of the first dynamical degree $d_1(f)$. The first dynamical degree is a bimeromorphic invariant and agrees with the spectral radius $\lambda(f)$ of the pullback $f^* \in \GL(H^2(X,\CC))$.
It is either equal to one or a so-called Salem number. The latter is a real algebraic integer $\lambda>1$ which is Galois conjugate to $\lambda^{-1}$ and such that all other conjugates lie on the unit circle.

Conjecturally, the smallest Salem number is Lehmer's number $\lambda_{10} \approx 1.17628$. The conjecture is true for the subset of Salem numbers coming from surface automorphisms: By a result of McMullen, the set of entropies of surface automorphisms has a spectral gap \cite[Theorem 1.2]{McMullen2007}, i.e.
\[h(f)=0 \quad \mbox{or} \quad h(f) \geq \log \lambda_{10}.\]  
Moreover, there exists a rational surface and a projective K3 surface with an automorphism of minimal entropy given by Lehmer's number $\lambda_{10}$ \cite{McMullen2016}.

The proof of existence for the K3 case works by constructing an integral Hodge structure $H$ of weight two and an effective Hodge isometry $f' \in \GL(H)$. Then the surjectivity of the period map and the strong Torelli theorem imply that there exist a K3 surface $X$, a marking $\eta\colon H^2(X,\ZZ) \to H$ and an automorphism $f\in \Aut(X)$ realising this data, i.e. with $\eta \circ f^* = f'\circ \eta$. Being of transcendental nature, the Torelli theorem is a mere existence result. It does not provide explicit equations.
In \cite{BrandhorstElkies23} the first named author and Elkies give explicit equations 
for an automorphism $f$ with minimal entropy of a K3 surface $X$. 
The approach there was to make use of the different elliptic fibrations of $X$.

Cantat \cite[Proposition 1]{Cantat1999} observed that if $d_1(f)>1$, then $X$ is bimeromorphic to either a rational surface, a 2-torus, a K3 surface or an Enriques surface. For $2$-tori the set of entropies is known thanks to results of Reschke \cite{Reschke2012,Reschke2017}. Oguiso and Yu \cite{OguisoYu20} show that for an Enriques surface $Y$ and an automorphism $\widetilde{f}\in \Aut(Y)$, its entropy $h(\widetilde{f})$ is either zero or at least $\log \tau_8\approx1.58234$ where $\tau_8$ is a root of the Salem polynomial
\[1 - x^2 - 2x^3 -x^4 +x^6 = 0.\]
Moreover, they prove the existence of a (complex) Enriques surface $Y_{OY}$ and an automorphism $\widetilde{f}_{OY} \in \Aut(Y_{OY})$ with entropy $h(\widetilde{f}_{OY}) = \log \tau_8$ by working on the Hodge-theoretic side. Again, this is an \textit{abstract} existence result, since the proof follows McMullen's strategy invokes the Torelli Theorem for Enriques surfaces.
Enriques surfaces are projective;
therefore Oguiso and Yu ask in \cite[Question 1.3]{OguisoYu20} for explicit equations of $Y_{OY}$ and $\widetilde{f}_{OY}$.
The purpose of this work is to answer their question and to popularise the methods developed along the way.

\medskip
The code for reproduction of the results in this article is available on \href{https://www.zenodo.org}{Zenodo} \cite{BrandhorstZach25code}. We will in the following
refer to the files and their names as uploaded there. 
\medskip

Recall that the universal cover of a complex Enriques surface $Y$ is a K3 surface $X$ and that the projection $X \to Y$ is of degree $2$. Let $\iota \in \Aut(X)$ denote the covering involution (in other words, the non-trivial deck transformation). It is fixed-point free.
Conversely, a fixed-point free involution on a K3 surface $X$ is called an Enriques involution, since the quotient $X/\iota$ is an Enriques surface. We denote by $X_{OY}$ the universal cover of $Y_{OY}$.
\begin{figure}
  \centering
  \includegraphics[scale=0.2]{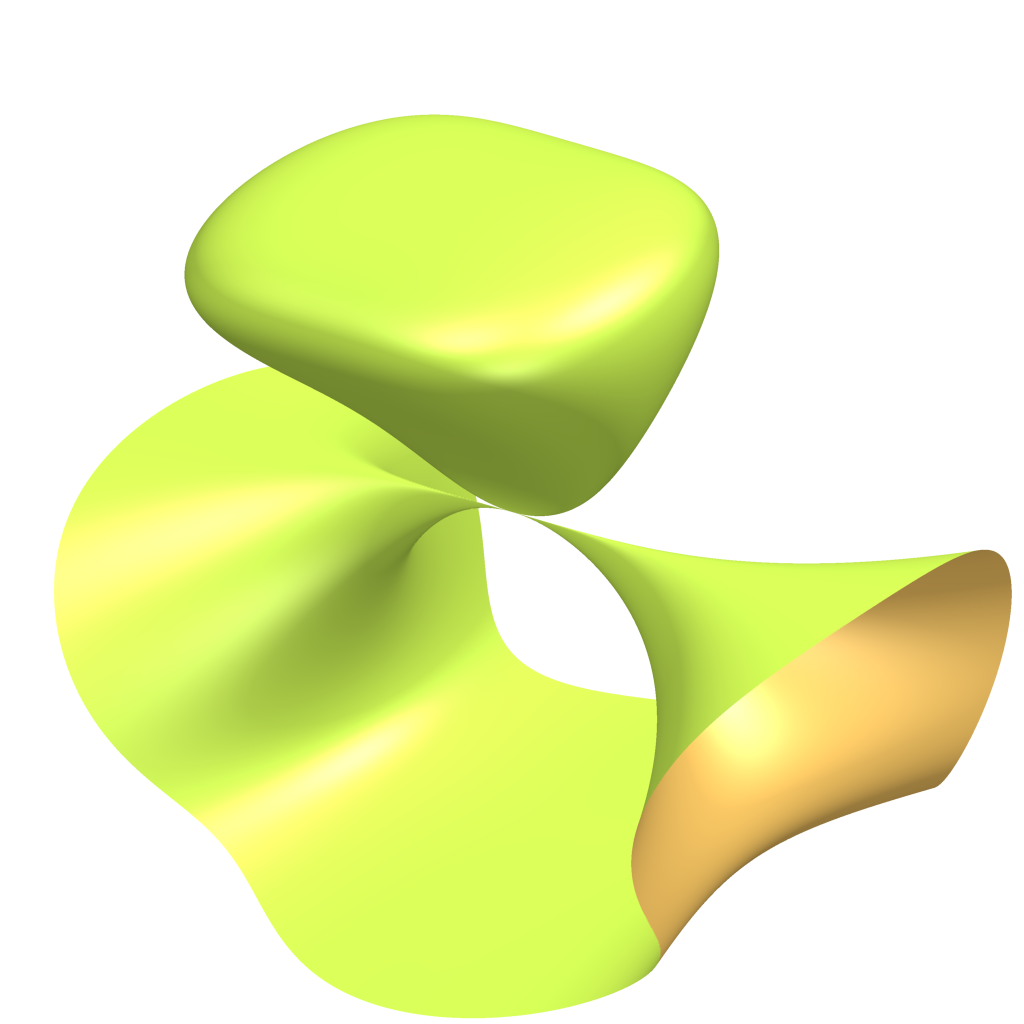}
  \label{fig:EnriquesSurfacePlot}
  \caption{The real points of the affine birational model for the Enriques surface sought for by Oguiso and Yu}
\end{figure}

\begin{theorem}\label{thm:main}
 The Oguiso-Yu Enriques surface is birational to the surface with equation
 \[\tilde y^2=s\tilde x^4+s^7-s^3.\]
 Its universal cover $X_{OY}$ is the elliptic K3 surface whose weierstrass equation is
 \[y^2 = x^3 + (1 - t^8)x.\]
 Translation by the  $2$-torsion section $(x,y) = (0,0)$ composed with $(x,y,t)\mapsto (x,y,-t)$ defines an Enriques involution $\iota_{OY}$. Then $Y_{OY} = X_{OY}/\iota_{OY}$ is the Oguiso-Yu Enriques surface.

 It admits an automorphism $\widetilde{f}_{OY
 } \in \Aut(Y_{OY})$ of entropy $h(\widetilde{f}_{OY
 }) = \log \tau_8$ which is defined over the number field $K=\QQ(\alpha)$ 
  where $\alpha$ is algebraic over $\QQ$ with minimum polynomial 
 \[a^8+6a^4+1=0.\]
 The automorphism is described in terms of
 its lift $f_{OY}$ to $X_{OY}$. 
 The pullback of the coordinate $t$ is given by
 \[
%
   f_{OY}^*(t) = \frac{x + (\frac{1}{2}\cdot \alpha^7 + 3\cdot \alpha^3 - \frac{1}{2}\cdot \alpha)\cdot y\cdot t + (\frac{1}{2}\cdot \alpha^6 + \frac{5}{2}\cdot \alpha^2)\cdot t^6 + (-\frac{1}{2}\cdot \alpha^6 - \frac{5}{2}\cdot \alpha^2)\cdot t^2}{x\cdot t^3 + (\frac{1}{2}\cdot \alpha^7 + 3\cdot \alpha^3 - \frac{1}{2}\cdot \alpha)\cdot y - t^5 + t}
\]
\end{theorem}

The equations of the pullback of $x$ and $y$ are far too large to display here.
They can be computed from the file \verb|reproduction_final_result.jl|.

\medskip
The benefit of concrete equations is that we can now actually map points, reduce the surface to positive characteristic and specify a field of definition of surface and automorphism.

\begin{corollary}
  The surface $Y_{OY}$ is defined over the field of rational numbers $\QQ$ and its automorphism of minimal entropy $\tilde f_{OY}$ is defined over the degree $8$ extension $\QQ(\alpha)$.
\end{corollary}

\begin{corollary}
 For all prime numbers $p>2$ there exists an Enriques surface $Y_p$ over
 $\FF_p$ and an automorphism $f_p \in \Aut(Y_p \times \Spec \FF_q)$ of entropy $\tau_8$ with $q=p^n$ and $n\leq 8$.
\end{corollary}
\begin{proof}
  Our model for $X_{OY}$ is defined over $\ZZ$ and has good reduction
  away from $2$. The Enriques involution is also defined over $\ZZ$. Since $p>2$, $\iota$ reduces to an
  Enriques involution $\iota_p$ on $X_p:= X_{OY}\times \Spec \FF_p$. The
  quotient $Y_p:=X_p/\iota_p$ is defined over $\FF_p$. Likewise $f_{OY}$
  is defined over the localization $\mathcal{O}_K[2^{-1}]$ and reduces modulo primes of $\mathcal{O}_K$ over $p$.
\end{proof}

\noindent
\textbf{Question:}  Is $\tau_8$ is the minimum entropy of an Enriques surface in characteristic $p$?

\subsection{Outline}
\footnote{A first draft of this outline was created using ChatGPT4o. Then it was partially rewritten by the authors.}
Section \ref{sec:preliminaries} briefly reviews K3 and Enriques surfaces, their Hodge structures, Torelli theorems, automorphism groups, and elliptic fibrations, as well as the role of lattice theory in classifying surface automorphisms.
The next Section \ref{sec:HodgeAndLatticeTheory} focuses on identifying
the Oguiso-Yu K3 surface $X_{OY}$
by searching a database for
K3 surfaces with a matching transcendental lattice.
The Enriques involution and automorphism group structure are analyzed,
confirming that the correct surface has been found.
Section \ref{sec:Construction} derives explicit equations by first representing $X_{OY}$  as a double cover $X_h \to \PP^1 \times \PP^1$, then identifying an elliptic fibration and computing a Weierstrass model $X_{f_1}$. The automorphism of minimal entropy is then constructed in Section \ref{sec:ConstructionAutomorphism} explicitly, and its action on the Néron-Severi lattice is verified.
In Section \ref{sec:Implementation}, we report on the framework for algebraic schemes used to carry out the computations and verify the end result.
Among others, the computational treatment of sheaves, birational maps, specialization, intersection theory and pushforward of divisors is explained in the context of a scheme represented by an atlas of charts.

\subsection*{Acknowledgements}
Special thanks goes to Matthias Schütt for discussions and suggesting the equation for the Enriques surface after reading a first version of this paper.
The first named author would like to thank Janko B\"ohm for sharing his insights.
We would further like to thank Anne Fr\"uhbis-Kr\"uger for collaboration and numerous
helpful discussions on the schemes package in \Oscar, as well as Max Horn, Claus Fieker, 
and Tommy Hofmann for their enduring support in integrating our work into the \Oscar 
codebase. Finally, thanks go to Hans Sch\"onemann for resolving various issues 
encountered with the Singular backend in \Oscar and again to Tommy Hofmann 
for providing numerous improvements on number-theoretic topics 
throughout the preparations of this work. 
This work was funded by the Deutsche Forschungsgemeinschaft (DFG, German Research Foundation) – Project-ID 286237555 – TRR 195.

\section{Preliminaries}\label{sec:preliminaries}
Since the universal cover of an Enriques surface is a K3 surface, we start with a brief review of K3 surfaces.
\subsection{K3 surfaces}
A \textit{complex K3 surface} is a compact, connected, complex manifold $X$ with vanishing irregularity and trivial canonical bundle
\[ h^1(X,\mathcal{O}_X)=0, \quad \omega_X\cong \mathcal{O}_X.\]

 The cup product gives the second singular cohomology group $H^2(X,\ZZ)$ the structure of a quadratic $\ZZ$-lattice, which turns out to be even, unimodular, and of signature $(3,19)$. Such a form is unique up to isometry.

Since K3 surfaces are simply connected, linear, algebraic, and numerical equivalence agree and $\Pic(X)=\NS(X)=\Num(X)$.
From the exponential sequence, we thus have an embedding
$c_1\colon \Pic(X) \to H^2(X,\ZZ)$ whose image is the subgroup of integral $(1,1)$-classes $H^{1,1}(X)\cap H^2(X,\ZZ)$ by Lefschetz' $(1,1)$-theorem.
The smallest primitive sublattice of $H^2(X,\ZZ)$ whose complexification contains the period $H^{2,0}(X)$, is called the \textit{transcendental lattice} $T(X)$. A K3 surface is not necessarily projective, but if it is, then $T(X)=\NS(X)^\perp$ in $H^2(X,\ZZ)$.

Let $X$ be a complex K3 surface.
The space of global $2$-forms $H^0(X,\omega_X)\cong H^{2,0}(X)$ is $1$-dimensional and generated by a nowhere vanishing symplectic form $\sigma_X$.
We denote by $\Aut_s(X)$ the kernel of the natural homomorphism
\[\Aut(X) \to \GL(\CC \sigma_X).\]
Its elements are called \textit{symplectic} automorphisms because they preserve the symplectic form $\sigma_X$. The other ones are called \textit{non-symplectic}. For $G\leq \Aut(X)$ set $G_s = \Aut_s(X)\cap G$. Note that if $X$ is projective, then $\Aut(X)/\Aut_s(X)$ is a finite cyclic group.

Complex K3 surfaces are governed to a large extent by their Hodge structure.
Let $X$ and $Y$ be two K3 surfaces. A \textit{Hodge-isometry} is a linear map
\[\phi\colon H^2(Y,\ZZ) \to H^2(X,\ZZ)\]
which preserves the intersection form and the Hodge structures, i.e. $\phi_\CC(H^{2,0}(Y))=H^{2,0}(X)$.
The weak Torelli theorem states that $X$ and $Y$ are isomorphic, if and only if $H^2(X,\ZZ)$ and $H^2(Y,\ZZ)$ are Hodge-isometric.

To formulate the strong Torelli theorem, we need one more definition.
The Hodge-isometry $\phi$ is called \textit{effective}, if it satisfies one of the following three equivalent conditions:
\begin{enumerate}
 \item $\phi$ maps effective classes on $Y$ to effective classes on $X$,
 \item $\phi$ maps K\"ahler classes to K\"ahler classes,
 \item $\phi$ maps a single K\"ahler class to a K\"ahler class.
\end{enumerate}
If $X$ is projective, then in (2) and (3) we may replace K\"ahler by ample.
\begin{theorem}[Strong Torelli theorem]
 Let $X,Y$ be K3 surfaces and
 \[\phi\colon H^2(Y,\ZZ) \to H^2(X,\ZZ)\]
 an effective Hodge isometry. Then there exists a unique isomorphism $f\colon X \to Y$ such that $\phi = f^*$.
\end{theorem}
Since the K3 surfaces covering Enriques surfaces are projective, we shall from now on assume that $X$ is projective and continue with the ample cone in place of the K\"ahler cone.
The connected component of 
\[
  \{x \in \NS(X)\otimes \RR \mid x^2>0\}
\]
which contains an ample class, is called the \textit{positive cone} $\mathcal{P}_X$. Let 
\[
  \Delta(X)=\{x \in \NS(X) \mid x^2=0\}.
\]
By the Riemann-Roch theorem we have that for $r \in \Delta(X)$ either $r$ or $-r$ is effective.
Set $\Delta^+(X) = \{r \in \Delta(X) \mid r \mbox{ is effective}\}$.
Then $h \in \mathcal P_X$ is ample if and only if $h.r>0$ for all $r \in \Delta^+(X)$.
More generally, the connected components of $\mathcal P_X\setminus \bigcup \{r^\perp \otimes \RR\mid r \in \Delta(X) \}$ are called Weyl-chambers. The ample cone is one such chamber. The Weyl group $W(X) \subseteq O(\NS(X))$ is generated by the reflections at the hyperplanes $r^\perp$. It acts simply-transitively on the set of Weyl-chambers.

Since $\Delta^+(X)$ is often infinite, it requires some work to control the ample cone. Useful tools are Vinberg's algorithm \cite{Vinberg72} and the algorithms described in \cite{Shimada24}.

\subsection{Enriques surfaces}
A \textit{complex Enriques surface} $Y$ is a compact complex manifold with $h^1(Y,\mathcal{O}_Y)=0$ and $b_2(Y)=10$. Its canonical bundle is $\omega_Y \neq \mathcal{O}_Y$, but $\omega_Y^{\otimes 2} \cong \mathcal{O}_Y$. Since the canonical bundle is $2$-torsion, it induces a $2:1$ cover $X \to Y$, which turns out to be the universal covering with $X$ a K3 surface. Since $h^{1,1}(Y)=b_2(Y)$, it follows from Lefschetz' Theorem on $(1,1)$-classes and the Kodaira Embedding theorem that
all Enriques surfaces are projective.

\subsection{Elliptic Fibrations}
An \textit{elliptic fibration} on a smooth projective surface $S$ is a morphism $S \to C$ to a curve $C$ whose general fiber is an elliptic curve. If the map admits a section, we call the fibration \textit{Jacobian}.

An elliptic fibration is called \textit{relatively minimal}, if its fibers do not contain any $(-1)$-curves. The singular fibers of a relatively minimal elliptic fibration have been classified by Kodaira into various types and they can be computed from a Weierstrass model by Tate's algorithm (see e.g. \cite{Silverman94}). The dual graph of a reducible fiber is an extended ADE-diagram.

For K3 and Enriques surfaces it is known that $C=\PP^1$.
Elliptic fibrations on Enriques surfaces are never Jacobian; on a K3 surface they may or may not be Jacobian. If it is, then its generic fiber admits a Weierstrass model
\[y^2 + a_1xy + a_3 y = x^3+a_4x^2+a_2x+a_6\]
with $a_i=a_i(t) \in \CC(t)$ of degree at most $2i$ (and further conditions assuring that it is smooth and minimal).

\subsection{Lattices}
Let $L$ be a quadratic integer lattice and $f \in O(L)$ an isometry.
The invariant and co-invariant lattices of $f$ are
\[L^f=\{x \in L \mid f(x)=x\} \qquad L_f:= (L^f)^\perp.\]
We follow the geometric convention that ADE-root lattices are negative definite.
\section{The Hodge- and lattice theory of the surface}
\label{sec:HodgeAndLatticeTheory}
In this section we deduce as much information as we can from the
Hodge-theoretic model of the sought for Enriques surface $Y_{OY}$ 
and its automorphism of minimal entropy, in order to guide our search for explicit 
equations.\\

\subsection{The construction of Oguiso and Yu}
We review the construction of Oguiso and Yu in \cite[Theorem 8.1 and proof]{OguisoYu20}.
Denote by $Y_{OY}$ the Oguiso-Yu Enriques surface, by $\widetilde{f}_{OY}$ an automorphism of entropy $h(\widetilde{f}_{OY})=\log \tau_8$,
by $\iota_{OY}$ the covering involution of the universal cover $X_{OY} \to Y_{OY}$
of $Y_{OY}$ by the K3 surface $X_{OY}$,
and by $f_{OY}\colon X_{OY}\to X_{OY}$ a lift of $\widetilde{f}_{OY}$.
Oguiso and Yu construct
\begin{enumerate}
 \item[a)] an even unimodular lattice $\Lambda_{K3}$ of signature $(3,19)$
 \item[b)] an isometry $i' \in  O(\Lambda_{K3})$ with  $i'^2=1$ and invariant lattice $\Lambda_{K3}^{i'} \cong U(2) \oplus E_8(2)$,
 \item[c)] an isometry $f' \in  O(\Lambda_{K3})$ with $i'\circ f' = f' \circ i'$ and spectral radius $\lambda(f')=\tau_8$.
 \end{enumerate}
Next, it is proven that there exists a Hodge structure of weight two on $\Lambda_{K3}$, preserved by both $i'$ and $f'$. Moreover, these two morphisms preserve the same Weyl-chamber 
and in particular, they are effective.
Then it follows from the surjectivity of the period map and the strong Torelli theorem that there exist
\begin{enumerate}
 \item[1)] a K3 surface $X_{OY}$,
 \item[2)] an Enriques involution $\iota_{OY} \in \Aut(X_{OY})$,
 \item[3)] an automorphism $f_{OY} \in \Aut(Y_{OY})$,
 \item[4)] a marking $\eta\colon H^2(X_{OY},\ZZ) \xrightarrow{\sim} \Lambda_{K3}$ such that $\eta \circ \iota_{OY}^* = i' \circ \eta$ and $\eta \circ f_{OY}^* = f' \circ \eta$.
\end{enumerate}
Finally, $f_{OY}$ descends to $Y_{OY} = X_{OY}/\iota_{OY}$ since it commutes with the Enriques involution $\iota_{OY}$.

\medskip
 All we get explicitly from this construction are a), b), and c) in the form of three $22\times 22$ integer matrices representing the gram matrix of $\Lambda_{K3}\cong \ZZ^{22}$ and the isometries $i'$ and $f'$.
 It is our task to find 1) to 4).

\subsection{The transcendental lattice}
\label{subsec:TranscendentalLattice}
By \cite[Theorem 8.1]{OguisoYu20} the transcendental lattice $T(X_{OY})$ is isometric to $L(4)$ with $L$ an odd unimodular lattice of signature $(2,2)$. Furthermore, $f_{OY}$ acts with order $8$ on $T(X_{OY})$.
Since $\varphi(8)=4$ (the Euler-phi-function), 
the characteristic polynomial of $f_{OY}|T(X_{OY})$ is $x^4+1$.
The complexification $T(X_{OY})\otimes \CC$ contains the period
\[H^{2,0}(X_{OY})\cong H^0(X_{OY},\Omega^2_{X_{OY}})= \CC \sigma_{OY},\] which is preserved by $f_{OY}$. Hence, it must be one of the four eigenspaces of $f_{OY}|T(X_{OY})\otimes \CC$. Moreover, we know that $\sigma_{OY}.\overline{\sigma_{OY}}>0$. This leaves us the possibility of only two complex conjugate eigenspaces which could contain $\sigma_{OY}$. Which one it is, 
is decided by the eigenvalue.
We see that the action of $f_{OY}$ on $H^2(X_{OY},\ZZ)$ and the eigenvalue of the eigenspace $H^{2,0}(X) \subseteq T(X_{OY})\otimes \CC$ determine the Hodge structure of $X_{OY}$. The next theorem shows that  under suitable assumptions the isomorphism class of the action on $T(X)$ and the eigenvalue determine the Hodge structure up to isomorphism.
\begin{theorem}\label{thm:uniqueX}\cite[Theorem
  1.2]{brandhorst:purely-nonsymplectic} Let $X_i$, $i=1,2$ be two complex
  K3 surfaces and $f_i \in \Aut(X_i)$ automorphisms with $f_i^* \sigma_i=
  \zeta_n \sigma_i$ where $\zeta_n = \exp(2\pi i/n)$, $n\in \NN$ and $\sigma_i
  \in H^0(X_i,\Omega_{X_i}^2)$ is a non-zero holomorphic $2$-form on $X_i$.
  Let $I_i$ be the kernel of the natural map given by \[\ZZ[\zeta_n]
  \to O(\NS(X_i)^\vee/\NS(X_i)), \quad \zeta_n \mapsto f_i^*.\] If $f_1$
  is of finite order, then $X_1 \cong X_2$ if and only if the ideals $I_1$
  and $I_2$ are equal.  
\end{theorem}

At this point it is plausible to search for some K3 surface with a non-symplectic automorphism of order $8$, having the same transcendental lattice as $X_{OY}$.

\medskip
A very general K3 surface with a purely non-symplectic automorphism of order $8$ has $2$-elementary transcendental lattice
and $\rank T(X) = 4 r$, a multiple of $4$.
Thus we have to look for a special member in a family of K3 surfaces with a non-symplectic automorphism of order $8$. One way to force extra algebraic cycles (hence forcing $r=1$) is to require the existence of additional symplectic automorphisms. We can therefore hope to find $X_{OY}$ among the set of K3 surfaces $X$ admitting a group $G\leq \Aut(X)$ with $G/G_s \cong C_8$ and its subgroup of symplectic automorphisms $G_s \neq 1$.

The database \cite{BrandhorstHofmann23} consists of a complete set of representatives of pairs $(X,G)$ where $X$ is a K3 surface and $G\leq \Aut(X)$ up to deformation and isomorphism. Each entry of the database corresponds to a moduli space of such pairs $(X,G)$. The moduli space is encoded by an even unimodular lattice $\Lambda_{K3}$ and a (satured and effective) group $G' \leq O(\Lambda_{K3})$, given in terms of matrices.
Via a marking (which we do not have explicit access to) $\Lambda_{K3}$ can be identified (non-uniquely!) with $H^2(X,\ZZ)$ and $G'$ with the image of $G \to O(H^2(X,\ZZ))$. For a generic member $X$ of the moduli space this allows to compute the isomorphism classes of $\NS(X)$, $T(X)$ and various further sublattices associated to the action of $G$ and $G_s$.

The database contains $108$ pairs $(\Lambda_{K3},G')$ with $G_s \neq 1$ and $|G/G_s|=8$. Among these, $66$ have $\rank T(X)=4$. They have $\det T(X) \in \{4,16,64,256,324\}$.
We are interested in $\det T(X)=256$. There are $18$ such database entries and all of them have $T(X) \cong T(X_{OY})$. They satisfy that $|G_s| \in \{4,8,16,32\}$. We will take $G_s$ as large as possible in the expectation that additional symmetry helps us to find the equations for $X$ and algebraic cycles in $\NS(X)$ and in the hope that $G$ contains the Enriques involution $\iota_{OY}$ we seek.

\begin{proposition}\label{proposition:IdentifyXG}
 Let $X$ be a complex K3 surface and $G \leq \Aut(X)$ such that $|G|=256$
 and $G/G_s \cong C_8$.
 Then the pair $(X,G)$ is unique up to isomorphism and given by the entry  \verb|40.8.0.1| in the database \cite{BrandhorstHofmann23}. Moreover, $G_s \cong (Q_8 * Q_8)$ and
 $X \cong X_{OY}$.
\end{proposition}

\begin{proof}
Only the entry numbered \verb|40.8.0.1| has a group action by a group $G$ of order $256$ and $G/G_s\cong C_8$. The corresponding moduli space consists of a single point and therefore the pair $(X,G)$ with the given conditions is unique up to isomorphism.
  Using the Hodge-theoretic model in the database, one computes that $T(X)$ is in the same genus as $T(X_{OY})$. The rescaled lattices $T(X)(1/4)$ and $T(X_{OY})(1/4)$ are both odd unimodular lattices of signature $(2,2)$. By the classification of indefinite odd unimodular lattices in terms of their signature (see e.g. \cite[Chapter V \S 2.2 Theorem 4]{Serre78}), they are abstractly isomorphic.

Let $gG_s$ be a generator of $G/G_s$.
We want to apply Theorem \ref{thm:uniqueX} with $X_1=X_h$, $X_2=X_{OY}$, $f_1=g^{k_1}$ and $f_2 = f_{OY}^{k_2}$ for $k_1$ and $k_2$ suitable odd integers such that $f_i$ acts by $\zeta_8$ on the space of holomorphic $2$-forms $H^{2,0}(X_i)$ for $i=1,2$.
Since $T(X)\cong T(X_{OY})$ as lattices,
the corresponding ideals $I_i \leq \ZZ[\zeta_8]$ are both of norm $N(I_i) = \det T(X)=4^4$. The prime $2$ is completely ramified in the extension $\QQ(\zeta_8)/\QQ$. Therefore, there is a unique ideal $I\leq \ZZ[\zeta_8]$ of norm $4^4=\det T(X)$. We see that $I_1=I_2$ and by Theorem \ref{thm:uniqueX} $X\cong X_{OY}$.
\end{proof}
Next, we need to search for clues where to find the Enriques involution $\iota_{OY}$ giving $Y_{OY}$.
By chance, it turns out that $\iota_{OY}$ is contained (up to conjugacy) in $G$.

\begin{remark}\label{remark:ohashi}
Let $X$ be a K3 surface and $i,j \in \Aut(X)$ be two Enriques involutions. The Enriques surfaces $X/i$ and $X/j$ are isomorphic if and only if $i$ and $j$ are $\Aut(X)$-conjugate \cite{ohashi}.
A necessary condition for this is that the algebraic coinvariant lattices $\NS(X)_i$ and $\NS(X)_j$
  are isomorphic. Note that $\NS(X)_i$ of $i$ does not contain any vectors $x$ with $x^2=-2$,
  since $i$ preserves the nef cone.
\end{remark}
By a direct computation with the Hodge-theoretic model provided by Oguiso and Yu,
one finds that $\NS(X_{OY})_{\iota_{OY}}$ is of determinant $2^{10}$ and lies in the genus $\mathrm{II}_{0,8}(2^{-2} 4^{-4}_4)$. An application of Kneser's neighbor method \cite[Kapitel IX]{Kneser02} (see \cite{Voight23} for a survey), as shipped with Oscar \cite{oscar}, shows that this genus contains a total of $3$ isometry classes, two of which have maximum $-2$ and one maximum $-4$. The unique one of maximum $-4$ must be the one isometric to $\NS(X_{OY})_{\iota_{OY}}$.
\begin{lemma}\label{lem:enriques_invol}
 Let $(X,G)$ be the symmetric K3 surface of id \verb|40.8.0.1|.
 Then $G$ contains up to conjugacy exactly two Enriques involutions $j_1$ and $j_2$.
 Their coinvariant lattices lie in the genera
 \[\NS(X)_{j_1} \in \mathrm{II}_{0,8}(2^{-2} 4^{-4}_4) \mbox{ and } \NS(X)_{j_2} \in \mathrm{II}_{0,8}(2^{-4} 4^{-4}_4).\]
The involution $j_1$ is contained in the center $C(G)\cong C_2 \times C_2$, but $j_2$ is not.
\end{lemma}
\begin{proof}
 An involution $j \in G$ is an Enriques involution if and only if
 its fixed lattice $F:=H^2(X,\ZZ)^j$ is isometric to $E_{10}(2)$ if and only if $F(1/2)$ is an even unimodular lattice (of signature $(1,9)$) \cite{Keum90}. Since the database contains abstractly the action of $G$ on $H^2(X,\ZZ)$, the conditions can be checked on a complete set of representatives of the conjugacy classes of involutions in $G$, yielding that $G$ contains two Enriques involutions, and that one of them is in the center. The coinvariant lattices and their genus can be computed explicitly.
 These computations are carried out in the file \verb|lattice_computations.jl|.
 \end{proof}

\begin{proposition}\label{prop:IdentifyInvolution}
 Let $X:=X_{OY}$ and $j \in \Aut(X)$ an Enriques involution with
 \[\NS(X)_{j} \in \mathrm{II}_{0,8}(2^{-2} 4^{-4}_4).\]
 Then $j$ is $\Aut(X)$-conjugate to $\iota_{OY}$.
\end{proposition}
\begin{proof}
 Using the theory of Nikulin and Miranda--Morrison \cite{nikulin, miranda_morrison:embeddings} and its \Oscar implementation \cite{BrandhorstHofmann23,brandhorst-veniani} one checks that any two primitive sublattices $F,F' \leq \NS(X)$ with $F \cong F' \cong \NS(X)_j$ are in the same $O(\NS(X))$-orbit.

Let $O(\NS(X),j)$ be the centraliser of $j$ in $O(\NS(X))$.
A direct (but non-trivial) computation based on Muller's Oscar package \cite{muller:quadform_and_isom, BrandhorstHofmann23} shows that the composition
\[O(\NS(X), j) \to O(\NS(X)) \to O(\NS(X)^\vee/\NS(X))\]
is surjective.
It follows from Shimada-Veniani \cite[Theorem 3.1.9]{shimada-veniani} that any Enriques involution $k \in \Aut(X)$ with
$\NS(X)_k \cong \NS(X)_j$ is $\Aut(X)$-conjugate to $j$. Take $k=\iota_{OY}$ to obtain the statement.
See the file \verb|lattice_computations.jl| for the computations.
\end{proof}

\begin{proposition}\label{prop:PairOfGroupAndSurfaceIsUnique}
Let $(X,G)$ be the symmetric K3 surface of id \verb|40.8.0.1| and denote by $j$ the unique Enriques involution in the center of $G$. Then
  $X/j \cong Y_{OY}=X_{OY}/\iota_{OY}$.
 \end{proposition}
\begin{proof}
By Lemma \ref{lem:enriques_invol} $\NS(X_{OY})_{\iota_{OY}}\cong \NS(X)_j$.
The statement of this proposition follows with Proposition \ref{prop:IdentifyInvolution} and Remark \ref{remark:ohashi}.
\end{proof}
To aid the search for a projective model of $X$ and $G$, we look at $G$-invariant ample divisors on $X$ and their properties.
Note that since $G$ is finite, there exists a $G$-invariant ample class. Since $\rank \NS(X)^G=1$, the invariant lattice is spanned by a unique primitive ample class $h \in \NS(X)^G$.
We have $h^2=4$. Furthermore, we know that $h$ is basepoint free,
because, by a direct calculation, the set
\[\{e \in \NS(X) \mid e^2=0, h.e=1\}= \emptyset\] 
is empty and hence $h$ induces a morphism $\phi: X \to \mathbb{P}^3$, cf. \cite{SaintDonat74}.
In general, the ample class $h$ can belong to one of the following three types according to 
\cite{SaintDonat74}:
\begin{enumerate}
 \item[A)] very ample: $\phi$ is an embedding
 \item[H)] hyperelliptic: $\phi$ is a $2:1$ cover onto a quadric hypersurface $Q\cong \PP^1\times \PP^1$ branched over a divisor of bidegree $(4,4)$.
 \item[U)] unigonal: the generic fiber of $\phi$ is an elliptic curve and its image is a twisted cubic.
\end{enumerate}

\begin{proposition}\label{prop:hyperelliptic}
 Let $(X,G)$ be the symmetric K3 surface of id \verb|40.8.0.1|.
 The (unique) $G$-invariant ample class $h \in \NS(X)^G$ is basepoint free hyperelliptic.
\end{proposition}
\begin{proof}
We have seen already that $h$ is basepoint free.
Now one computes that \[\{e \in \NS(X) \mid e^2=0, h.e=2\}= \{e_1,e_2\}.\]
By \cite{SaintDonat74} the linear system $|h|$ is hyperelliptic.
Note that we can write $h = e_1 + e_2$ for two nef isotropic classes $e_1$ and $e_2$ with $e_1.e_2 = 2$.
The two isotropic classes give the two elliptic fibrations $X \to \PP^1 \times \PP^1 \to \PP^1$ with the first map induced by $|h|$ and the second map the projection to one of the factors.
See the file \verb|lattice_computations.jl| for the computations.
\end{proof}
\begin{remark}
 The ample class $h \in \NS(X)^G$ is invariant under the Enriques involution $j \in G$. Hence $h$ descends to an ample class
 $h_Y$ of degree $2$ on $Y_{OY}\cong X/j$.
 The elliptic divisors $e_1,e_2 \in \NS(X)$ with $e_1+e_2=h$ induce elliptic fibrations on $X$. Since they are invariant under $j$, they descend to $Y$ and induce elliptic fibrations on $Y$.
  The morphism corresponding to the sum of the two half fibers $\pi_*(e_i)$ gives the so called Horikawa model of non-special Enriques surfaces.
  Their geometry is described for instance in \cite[VIII Proposition 18.1]{BarthHulekPetersVen04}.
\end{remark}
  Let us compute the type of the elliptic fibrations induced by $e_1$ and $e_2$.
\begin{proposition}\label{prop:elliptic_fibration_on_X}
Let $i \in \{1,2\}$. The elliptic fibration induced by $e_i$ on $X$ has $8$ reducible fibers of type $A_1$ and Mordell-Weil group
$\ZZ/2\ZZ \oplus \ZZ^8$. The stabiliser of $e_i$ in $G$ is of order $128$ and its orbit is $\{e_1,e_2\}$.
The elliptic fibration induced by $e_i$ on $Y= X/j$ has $4$ simple $A_1$ fibers.
\end{proposition}
\begin{proof}
 The description of the stabilisers and the orbit follows from the proof of Proposition 
  \ref{prop:hyperelliptic} as follows: Since $h$ is invariant under $G$, so is the set
 \[\{e \in \NS(X) \mid e^2=0, h.e=2\}=\{e_1,e_2\}.\]
 Since the invariant lattice is of rank $1$, $G$ must exchange $e_1$ and $e_2$ and the stabiliser of $e_i$ is of order $G/2=128$.
 Recall from \cite{shioda} that the trivial lattice $\Triv(e)$ of an elliptic fibration with fiber $e \in \NS(X)$ is the sublattice of $\NS(X)$ spanned by the zero section and fiber components. The Mordell-Weil group is isomorphic to $\NS(X)/\Triv(e)$.
 To compute the trivial lattice, one chooses any $s \in \NS(X)$ with $s^2=-2$ and $s.f=1$. Let $R$ be the root sublattice of $\langle e, s \rangle^\perp \subseteq \NS(X)$. Then the trivial lattice is isomorphic to $\ZZ \langle s, e \rangle \oplus R$. The irreducible root sublattices of $R$ give the ADE types of the singular fibers.
 For how to compute the fiber types on the Enriques surface $Y$ from the class of $e_i$ see \cite[Prop. 5.8]{brandhorst-gonzalez}.
\end{proof}

\begin{remark}
With the methods described in \cite[Section 5]{brandhorst-gonzalez} and the Hodge and lattice theoretic model of $(X,G)$, one computes the following:
Let $Y_{OY}$ be the Oguiso-Yu Enriques surface.
Modulo $\Aut(Y_{OY})$ it contains exactly $12$ elliptic fibrations.
The following table contains the ADE types of the reducible fibers and whether they are double, the
count of $\Aut(Y_{OY})$ orbits with the same fiber types, and the ramification index of the fibration in
the forgetful map $\mathcal{P}^0\colon \mathcal{M}_{En,0} \to \mathcal{M}_{En}$ of the moduli space of elliptic Enriques surface to the moduli space of Enriques surfaces.
\begin{center}
\begin{tabular}{cccc}
simple fibers & double fibers & count & ramification index\\
\hline
$1 A_1 + A_3$ & - & 1 & 128\\
$2 A_1 + A_3$ & - & 1 & 128\\
$3 A_1 + A_3$ & -  & 2 & 64\\
$6 A_1$ & -  & 1 & 48\\
$4 A_1$ & - & 1 & 48\\
$4 A_1$ & $2A_1$& 1 & 24\\
$ 2A_1+ D_4$ & - & 2 & 8\\
$ 2A_1+ D_4$ & $A_1$ & 1 & 4\\
$ 2A_1$ & $A_3+A_1$ & 1 & 2\\
$ 4A_1$ & $A_3$ & 1 & 1\\
\end{tabular}
\end{center}
As a sanity check we calculate that the weighted count of isomorphism classes of elliptic fibrations equals
\[128+128+2 \cdot 64+48+48+24+2 \cdot 8+4+2+1 = 527,\] which agrees with the degree of the forgetful map $\mathcal{P}^0$ -- as it should.
\end{remark}
Besides a projective model for $X_{OY}$ and an elliptic fibration, we also seek explicit generators for $\NS(X_{OY})$.
A common approach is to search for smooth rational curves of small degree. The lattice side will tell us up to which degree we have to search:
\begin{proposition}\label{prop:rational_curves_degrees}
   Let $(X,G)$ be the symmetric K3 surface of id \verb|40.8.0.1| and $h$ be the unique $G$-invariant ample class of $X$. Denote by $D_i$ the set of smooth rational curves $c$ on $X$ of degree $h.c = i$. Then $|D_1|=32$, $\rank \langle D_1 \rangle = 14$ and $D_1$ forms a single $G$-orbit.
The set $D_2$ is of cardinality $160$ and splits into three $G$-orbits of sizes $32$, $64$, and $64$.
Altogether, $D_1$ and the orbit of size $32$ in $D_2$ span a sublattice $L$ of $\NS(X)$ of rank $18$, determinant $2^2 4^4$ and index $[\NS(X):L] = 2$.
\end{proposition}
\begin{proof}
The classes in $\NS(X)$ of the smooth rational curves of $X$ are the facets of the nef cone, which itself is a fundamental chamber for the action of the Weyl group generated by the $(-2)$-hyperplanes. We use Vinberg's algorithm\cite{Vinberg72} to compute the facets of smallest distance to $h$ in hyperbolic space:
\[D_1 = \{x \in \NS(X) \mid h.x=1, x^2=-2\},\]
\[D_2 = \{y \in \NS(X) \mid h.y=2, y^2=-2, \forall x \in D_1: x.y\geq 0\}.\]
These sets can be computed by solving a suitable close vector problem.
See the file \verb|lattice_computations.jl| for the computations.
\end{proof}

At this point we are given on the one hand
$(\Lambda_{K3}, \iota', f')$ corresponding to $X_{OY}$, $\iota_{OY}$ and $f_{OY}$.
On the other hand from the database entry \verb|40.8.0.1| we get $L_{K3},G',j', h'$ corresponding to $X,G$, and $h$. Both sides come with implicitly defined Hodge-structures and Weyl-chambers.
Since $X_{OY}\cong X$, we need to reconcile these two points of views, i.e. we need to find an explicit effective Hodge isometry $\Lambda_{K3} \to L_{K3}$ that conjugates $\iota'$ to $j'$.
Ideally, it should map $f'$ to a matrix with small coefficients in our choice of basis for $\NS(X_{OY})\otimes \QQ$ from $D_1 \cup D_2$ (or rather $f'(h).h$ should be small).

Given the existing computational infrastructure it is easier to compute generators for the image of $\Aut(Y) \to \Num(Y)$ (where $Y=X/j$ and $\Num(Y)(2)$ is identified with $L_{K3}^j$) using Borcherds' method as described in \cite{brandhorst_shimada:tautaubar} and search for an automorphism with entropy $\log \tau_8$. This has the additional merit of possibly factoring $f_{OY}$ into isometries of zero entropy and hence a simpler geometry.

In the following proposition we collect our findings:
\begin{proposition}\label{prop:factor_OY}
  Let $Y=X/j$ be the Oguiso-Yu Enriques surface where $(X, G)$ is the entry \verb|40.8.0.1|
  from \cite{BrandhorstHofmann23}.
  Then $X$ carries an automorphism $f$ of entropy $\log \tau_8$ and commuting with $j$. It factors as
$f = l \circ g$ where $g \in G$ and $l$ preserves an elliptic fibration $\phi_l \colon X \to \mathbb{P}^1$ with
two $\widetilde{A}_1$ fibers and two $\widetilde{A}_3$ fibers.
The Mordell-Weil group of $\phi_l$ is $\ZZ^8$ and this fibration descends to $Y$ with a reducible fiber of type $A_1$ and a reducible fiber of type $A_3$.
\end{proposition}
We have collected enough information to guide the search for equations of the symmetric K3 surface $(X,G)$, for the elliptic fibration on $X$ given by the isotropic class $e_1$ and for generators of its Néron-Severi lattice $\NS(X)$. Once we have these, it remains to connect them by an effective Hodge isometry to the abstract K3 surface found in the database.

%
%
%
%


\section{Construction of the Enriques surface}\label{sec:Construction}

\subsection{The first model}
Recall from Proposition \ref{prop:PairOfGroupAndSurfaceIsUnique} that we have identified the universal 
cover $X_{OY} \to Y_{OY}$ of the sought for Enriques surface $Y_{OY}$ 
to be the entry $(X, G)$ with number \verb|40.8.0.1| from
\cite{BrandhorstHofmann23}.
In order to construct a first concrete model $X_h$ for $X_{OY}$, we start from the
hyperelliptic system $|h|$ from Proposition \ref{prop:hyperelliptic} with associated morphism:
\begin{equation}
  \label{eqn:SaintDonatSystem}
  \Phi \colon X_h \overset{2:1}{\twoheadrightarrow} Q \cong
  \PP^1 \times \PP^1 \hookrightarrow \PP^3
\end{equation}
which is a $2:1$ cover onto a quadric hypersurface $Q \cong \mathbb{P}^1 \times \mathbb{P}^1$,
branched over a divisor $B$ of bidegree $(4,4)$.
Since $|h|$ is ample, $\Phi$ does not contract curves and hence $B$ is smooth.

From the discussion in \cite[VIII \S 18]{BarthHulekPetersVen04} we know that the polynomial defining $B$ is anti-invariant under the Enriques involution.
Since $h$ is $G$-invariant, we know that $G$ acts on the linear system $|h|\cong \mathbb{P}(H^0(X,\mathcal{O}_X(h))^\vee)$ by projective linear transformations. The action of $G$ on $|h|$ has a kernel of order $2$ consisting of the covering involution of the map $\Phi$.
Given this context, it is natural to look for projective representations
of $G = (Q_8 * Q_8).C_8$ on $\PP^3$ and their relative invariants, i.e. polynomials preserved by $G$ up to scaling.
These can be enumerated using \Oscar.

Eventually, we find that there is a projective representation 
$\rho \colon G \to \PP GL(3; \QQ(\zeta_8))$ 
over the cyclotomic field generated by a primitive $8$-th root of unity $\zeta_8$. The representation $\rho$ is not faithful and its kernel is generated by an involution $\eta\in G$.
The image $\tilde G \subset \PP GL(3; \QQ(\zeta_8))$ 
of this representation is the projective matrix group generated by\footnote{Note that in accordance with the \Oscar convention matrices act from the right throughout this paper and a ``vector'' is a row-vector unless otherwise specified.}
\begin{equation}
  \label{eqn:GeneratorsGTilde}
  \widetilde{g}_1 =
  \begin{pmatrix}
    0 & 1 & 0 & 0\\
    -1 & 0 & 0 & 0\\
    0 & 0 & 1 & 0\\
    0 & 0 & 0 & -1
  \end{pmatrix}, \quad
  \widetilde{g}_2 =
  \begin{pmatrix}
    0 & 0 & 1 & 0 \\
    0 & 0 & 0 & 1 \\
    \zeta_8 & 0 & 0 & 0 \\
    0 & -\zeta_8^3 & 0 & 0 
  \end{pmatrix}.
\end{equation}
The center of $\widetilde{G}$ is of order $2$ and generated by $\widetilde{i}:=\widetilde{g}_1^2$.

The \textit{relative invariants} of $\rho$ are given by two polynomials 
\[
  q(x, y, z, w) = xy + zw, \quad b(x, y, z, w) = x^4 + y^4 - i z^4 - i w^4.
\]
The expectation now is that $q$ defines the quadric $Q$ from the system
(\ref{eqn:SaintDonatSystem}), 
embedded as $\PP^1 \times \PP^1 \hookrightarrow \PP^3$ via 
\[
  [(s:t), (u:v)] \mapsto (x:y:z:w) = (su: tv:tu:-sv)
\]
and that the vanishing locus of $b$ is the 
branch locus $B$ of $\Phi$.
We thus construct a $2:1$ covering from $b$ by 
setting 
\begin{equation}
  \label{eqn:LocalDefiningEquationOfX}
  \lambda^2 = b(1, y, z, w)
\end{equation}
for the affine coordinate $\lambda$ of the bundle $\PP(\OO \oplus \OO(-2))$ 
over $Q$.
Extending this in the obvious way to all other affine charts, we obtain a surface 
\[
  \xymatrix{
    X_h \ar@{^{(}->}[r] &
    \PP (\OO \oplus \OO(-2)) \ar[d] & 
    \\
    & Q \ar[dl]_{\pi_1} \ar[dr]^{\pi_2}& \\
    \PP^1 & & \PP^1
  }
\]
with two natural projections to $\PP^1$ induced from the rulings of 
$Q \cong \PP^1 \times \PP^1$. 
This construction is carried out in \verb|minimum_entropy_surface.jl|, lines 1--72.
By direct calculation we verify that the curve $B$ is smooth and, hence, 
by \cite[Chapter V Section 23]{BarthHulekPetersVen04}, $X_h$ is a K3 surface.

Let $\widetilde{j} \in C(\widetilde{G})$ be the involution in the center of $\widetilde{G}$. It is given by $\widetilde{j}(x,y,z,w) = (-x,-y,z,w)$.
It has $4$ isolated fixed points and none of them lies on the branch curve $B$. By \cite[Chapter V Section 23]{BarthHulekPetersVen04} the involution $\widetilde{j}$ lifts to a fixed-point-free involution $j \in G$ of $X_h$ and $Y:=X_h/j$ is the Horikawa model of a so called non-special Enriques surface. Further, the two rulings of $Q$ define elliptic fibrations $X_h \to Q \to \mathbb{P}^1$.

It is easy to see that both $\pi_1$ and $\pi_2$
induce elliptic-fibrations on $X_h$: The fiber $\pi_1^{-1}(\{p\}) = \{p\}\times \PP^1$
of any point $p \in \PP^1$ is a projective line. The preimage of that space 
in $X_h$ is a curve $C$ with a generically $2:1$-projection $C \to \PP^1$,
ramified over the points where $b = 0$. For a generic choice of $p$ the curve 
$C$ is smooth and there are four simple ramification points. Then the Hurwitz formula
for morphisms of smooth curves asserts that in fact $g(C) = 1$. In the following we will denote the composed 
maps $X_h \to \PP^1$ by $\pi_1$ and $\pi_2$, as well.

\medskip

Following Proposition \ref{proposition:IdentifyXG} we know that we have found the pair $(X,G)$ of id \verb|40.8.0.1| as soon as we show that $G/G_s \cong C_8$, i.e. that $G$ acts with order $8$ on $H^{2, 0}(X_h, \CC) = H^0(X_h, \Omega_{X_h}^2)$.
We quickly verify this via explicit calculation as follows.

\begin{lemma}
  A global non-vanishing holomorphic $2$-form $\sigma_{X_h} \in H^0(X_h, \Omega^2_{X_h})$ is given by the natural
  extension of 
  \[
    \frac{1}{2\lambda} \D t\wedge \D v
    =\frac{1}{4(v^4 - i)t^3} \D \lambda \wedge \D v
    =\frac{1}{4(t^4 - i)v^3} \D t \wedge \D \lambda
  \]
  in the defining chart of (\ref{eqn:LocalDefiningEquationOfX}) and the affine coordinates 
  $(t, v)$ of $\PP^1 \times \PP^1$.
\end{lemma}

\begin{proof} This is a slightly tedious, but straightforward computation:
%
%
The automorphism on $\PP^1 \times \PP^1$ induced by $\tilde g_2$ is given by 
\[
  \tilde g_2 \colon (t, v) \mapsto \left(\frac{1}{\zeta_8 t}, \zeta_8^3 v\right)
\]
in the standard affine coordinates $(t, v)$ of the chart $\{s \neq 0\} \cap \{u \neq 0\}$ of $\PP^1\times \PP^1$.
The defining equation (\ref{eqn:LocalDefiningEquationOfX}) pulls back to 
\begin{equation}\label{eqn:chart_tv}
  \lambda^2 = 1 + t^4 v^4 - i t^4 - i v^4
\end{equation}
in these coordinates and we can lift $\tilde g_2$ to an automorphism of $X_h$ by setting
\begin{equation}
  \label{eqn:LiftingG2}
  \tilde g_2' \colon (\lambda, t, v) \mapsto \left(\frac{\zeta_8\lambda}{t^2}, \frac{1}{\zeta_8 t}, \zeta_8^3 v\right).
\end{equation}
Indeed, for these choices we get 
\begin{eqnarray*}
  \left(\tilde g_2'\right)^*(1 + t^4 v^4 - i t^4 - i v^4) & = & 
  \left(1 + \frac{v^4}{t^4} + i\frac{1}{t^4} + i v^4\right) \\
  & = & 
  i \frac{1}{t^4} \left(-i t^4 - iv^4 + 1 + t^4v^4\right)\\
  & = & 
  i \frac{1}{t^4} \lambda^2 = \left(\frac{\zeta_8 \lambda}{t^2}\right)^2 = \left(\tilde g_2'\right)^* (\lambda^2).
\end{eqnarray*}
But then 
\[
  \left(\tilde g_2'\right)^* \sigma_{X_h} =
  \frac{t^2}{2\zeta_8 \lambda} \cdot \left(\frac{-1}{\zeta_8 t^2} \D t\right) \wedge \left(\zeta_8^{3} \D v\right)
  = -\zeta_8 \cdot \sigma_{X_h}
\]
so that the lift of $\tilde g_2$ clearly acts with order $8$ on the subspace generated by $\sigma_{X_h}$.

Alternatively, one may argue as follows: By direct calculation $(\tilde g_2')^4(\lambda,t,v) = (\lambda,t,-v)$.
Hence it has a fixed curve given by $v=0$ and must be non-symplectic, because a finite order symplectic automorphism has only isolated fixed points.
\end{proof}
It then follows from Proposition \ref{proposition:IdentifyXG} that $X_{OY} \cong X \cong X_h$ and from Proposition \ref{prop:PairOfGroupAndSurfaceIsUnique} that $X_h/j=:Y$ is isomorphic to the Oguiso-Yu
Enriques surface $Y_{OY}$ with an automorphism of minimum entropy.
\subsection{Search for sections}
\label{sec:FindingSections}
So far we have found the surface $X_h \cong X_{OY}$, the Enriques involution $\iota_{OY}$, and hence $Y_{OY} = X_{OY}/\iota_{OY}$. To hunt for equations of the automorphism $f_{OY}$ of minimum entropy, we need to find generators of the Néron-Severi lattice $\NS(X_h)$ and
a Weierstrass model $X_{e_1}$ of the elliptic fibration $\pi_1\colon X_h \to \mathbb{P}^1$; in particular a section.

Proposition \ref{prop:rational_curves_degrees} tells us that to generate a finite
index sublattice of $\NS(X_h)$, we have to search for smooth rational curves of degrees
$1$ and $2$ contained in $X_h$. Since $h = e_1+e_2$ (with $e_i$ the class
of a fiber of either one of the two elliptic-fibrations $\pi_i$) the
smooth \textit{rational} curves $C$ of degree $1=[C].h$ are contained in singular fibers of $\pi_1$ or $\pi_2$, hence easy to write down.
Thus it remains to search for rational curves $D \subseteq X_h$ of degree
$2=[D].h= [D].(e_1+e_2)$. There are two such cases which can occur.
Either $[D].e_1=0$ and $[D].e_2=2$ or vice versa and hence $D$ is contained
in one of the fibers. Or --and this is the other case-- $[D].e_1=[D].e_2=1$,
i.e. $D$ is a section of both fibrations.

In order to find sections $D$ of the maps $\pi_i \colon X_h \to \PP^1$ as divisors on $X_h$, we
reason as follows. Write 
$e_i = \Phi^*\epsilon_i$ for some divisors $\epsilon_i$ on $Q$. Then $1 = e_i.D = \Phi^*\epsilon_i.D = \epsilon_i.\Phi_*D$ for $i=1,2$ shows that $\Phi_*D \subset \PP^1 \times \PP^1$
is the graph of an isomorphism $\mathbb{P}^1 \to \mathbb{P}^1$.
Further $\Phi^*\Phi_*D\neq D$, since otherwise 
\[
  e_i.D = (\Phi^*\epsilon_i).(\Phi^*\Phi_*D)
  = \deg \Phi \cdot (\epsilon_i.(\Phi_*D))
  = 2 \cdot (\epsilon_i.(\Phi_*D))
\]
would be even. Hence $\Phi^*\Phi_*D$ must split into two components $D_1 + D_2$. This leads us to the following Ansatz. Let
\[
  \Delta_a = \{((s:t), (s:a\cdot t) \in \PP^1 \times \PP^1 \mid (s:t) \in \PP^1\},
\]
be a ``tilted diagonal''. 
If we restrict the defining equation $b$ of the branch locus to $\Delta_a$, 
then it might split as a square, depending on $a$
\[
  b(x, y, z, w)|_{\Delta_a} = b(s^2, a\cdot t^2, st, -a^2 \cdot st) = 
  s^8 + a^4 \cdot t^8 - i s^4 t^4 - ia^4 s^4 t^4 
  = s^8 - i\left(a^4 + 1\right) \cdot s^4 t^4 + a^4 \cdot t^8.
\]
This is a homogeneous polynomial of degree two in the variables $(s^4, t^4)$.
Hence, for either one of the eight solutions of the discriminant
\begin{equation}\label{eqn:MinPolyFieldExtension}
  (a^4 + 1)^2 + 4a^4 = a^8 + 6a^4 + 1 = 0
\end{equation}
we find two components $\{C_{j, 1}, C_{j, 2}\}_{j=1}^{8}$ 
of $\Phi^{-1}(\Delta_a)$. These are naturally sections of either one of the two fibrations $\pi_i$. 
In order to gain computational access to such solutions for $a$, we extend the
field 
over which we are working from $\QQ(\zeta_{8})$ to the splitting field $K$ of the above discriminant. 
The degree of the extension is $[K:\QQ(\zeta_8)]=2$.

\medskip
\noindent
These computations are carried out in lines 75-97 of \verb|minimum_entropy_surface.jl|. Note that in that script
we start out with the field $K$ from the beginning in order to avoid complicated base changes later on.

\begin{remark}\label{rem:FurtherComputationsInFirstConstruction}
  In this first construction of $X_h \cong X_{OY}$, we were also able to extract the components of the reducible fibers and realize the action by the group $G$ on $X_h$. As predicted by Proposition \ref{prop:rational_curves_degrees} the sublattice $L$ of $\NS(X_h)$ spanned by the curves of degree $2$ we found has rank $18$ and determinant $1024= [ \NS(X_h):L]\det(\NS(X_h))$.
\end{remark}

\subsection{Extracting a Weierstrass model}
The construction of $X_h \overset{2:1}{\longrightarrow} \PP^1 \times \PP^1 \to \PP^1$ as
a $2:1$-cover already suggests an entry point for the computation of a Weierstrass equation.
Passing to the affine chart $U$ where $x=1$ of the ambient space $\PP(\OO \oplus \OO(-2)|_Q)$ of $X_h$
with coordinates $(\lambda, y, z, w)$, we can use the defining equation $q=xy+zw$ for $Q$ to eliminate $y$. Then we are left with
\[
  \lambda^2 = 1 + z^4 w^4 - i z^4 - i w^4,
\]
which is the pullback of Equation \ref{eqn:chart_tv} to the respective affine chart over $\PP^1\times \PP^1$
discussed earlier. 
This is already an equation of a suitable format and classical formulas apply to transform it to a Weierstrass equation; see e.g. the work by A. Kumar \cite{Kumar14} for such procedures.
To this end we need to specify a
section, i.e. a $K(t)$-rational point in the associated elliptic curve. We choose this section
to be one of those constructed earlier.

The half-automated extraction of a globally minimal Weierstrass equation for an elliptic curve over $K(t)$, its simplification and the rational maps for the various identifications are carried out in lines 100-200 of
\verb|minimum_entropy_surface.jl|.
We end up with the Weierstrass equation 
\begin{equation}  \label{eqn:FirstWeierstrassEquation}
y^2 = x^3 - (t^8+1)x.
\end{equation}
\begin{remark}
  In hindsight, one could have arrived at this equation starting from
  Proposition \ref{prop:elliptic_fibration_on_X} which gives the existence
  of an elliptic fibration on $X_{OY}$ with $8$ fibers of type $A_1$
  and a $2$-torsion section.
  Since there is a non-symplectic automorphism $\alpha$ of order $8$
  preserving the fibers and the section, this in fact determines the
  Weierstrass equation uniquely: Looking at the action of $\alpha$ on
  $\NS(X_{OY})$ one sees that $\alpha$ acts with order $8$ on the the
  base $\mathbb{P}^1$. Up to a change of coordinates, $\alpha$ thus
  acts by $\zeta_8$ on $t$. Since there is a $2$-torsion section the
  Weierstrass equation is of the form $y^2=x^3+c_4(t)x$ with $c_4 \in
  K[t]$ of degree $8$, with $8$ simple roots giving the $A_1$-singular
  fibers and invariant under $\alpha$. Hence $c_4(t) = at^8+b$. After
  a change of coordinates one may assume $a=b=-1$.
\end{remark}
\medskip 
In order to be able to work with smaller extensions of $\QQ$,
we pass from the equation (\ref{eqn:FirstWeierstrassEquation})
to the elliptic surface with generic fiber $E_1$ defined by
\begin{equation}
  \label{eqn:ModifiedFirstWeierstrassEquation} 
  y^2 = x^3 + (1 - t^8)x,
\end{equation} 
which is isomorphic to the previous one, as they
are related by a change of Weierstrass coordinates
\[
  (x, y, t) \mapsto (-\zeta_{16}^4 \cdot x,\; -\zeta_{16}^2\cdot y,\; \zeta_{16} \cdot t)
\]
where $\zeta_{16}$ is a primitive $16$-th root of unity.
The associated elliptically fibered K3 surface $X_{e_1}$ is defined
over $\QQ$ and its trivial lattice over $\QQ(\zeta_8)$. However, in order to
realize all algebraic cycles for the N\'eron-Severi lattice $\NS(X_{e_1})$, we still need the extension $K/\QQ(\zeta_8)$ of degree two,
generated by the algebraic element $\alpha$
with minimal polynomial over $\QQ$ given in Equation \ref{eqn:MinPolyFieldExtension}.\\

\begin{proposition}\label{prop:EnriquesProjection}
The surface
\begin{equation}\label{eqn:YOY-equation}
 \tilde{y}^2=s\tilde{x}^4+s^7-s^3
\end{equation}
is birational to the Oguiso-Yu Enriques surface $Y_{OY}$.
\end{proposition}
Before we go into the proof, we review how to obtain an elliptic fibration of the covering K3 surface greater generality:
Let $\widetilde{\pi}\colon Y\to \mathbb{P}^1$ an elliptic fibration of an Enriques surface with double fibers over $P_1,P_2 \in \mathbb{P}^1$, and let $\kappa \colon X \to Y$ be the universal cover. Suppose that $\widetilde{\pi}$ admits a rational bisection.
According to \cite[Proposition 4.10.10]{cdl:EnriquesI} the Stein-factorization of
$\widetilde{\pi}\circ \kappa \colon X \to Y \to \mathbb{P}^1$  is given by an elliptic fibration
$\pi\colon X_{OY} \to \mathbb{P}^1$ followed by a $2:1$ cover $\PP^1 \to \PP^1$ ramified over $P_1,P_2$.

\begin{proof}[Proof of Proposition \ref{prop:EnriquesProjection}]
  We set
  $t^2=s$, $\tilde x=t x_1$, $\tilde{y}=t^3y_1$ in equation (\ref{eqn:YOY-equation}) 
  and obtain $y_1^2=x_1^4+t^8-1$. Using the rational point
  $(x_1, y_1) = (1, t^4)$
  we can transform the equation to a short Weierstrass form
  which one checks to be isomorphic to the one of eq. (\ref{eqn:ModifiedFirstWeierstrassEquation}).
  Putting everything together, we obtain obtain the universal covering morphism $\kappa_{OY} \colon X_{OY}\to Y_{OY}$ by 
  means of the following rational expressions in the Weierstrass coordinates of $X_{OY}$:
  \[
    \begin{pmatrix}
      \tilde x \\
      \tilde y \\
      s
    \end{pmatrix}
    \mapsto
    \begingroup
    \renewcommand*{\arraystretch}{2}
    \begin{pmatrix}
      \frac{-x\cdot t^9 + x\cdot t + \frac{1}{2}\cdot \sqrt 2\cdot y\cdot t^5 - t^9 + t}{x + \frac{1}{2}\cdot \sqrt 2\cdot y\cdot t^4 + t^8 - 1}\\
      \frac{-x^3\cdot t^{15} - 3\cdot x^2\cdot t^{15} + 3\cdot x^2\cdot t^7 + \frac{1}{2}\cdot y^2\cdot t^{15} + \sqrt 2\cdot y\cdot t^{19} - 3\cdot \sqrt 2\cdot y\cdot t^{11} + 2\cdot \sqrt 2\cdot y\cdot t^3 + t^{23} - 2\cdot t^{15} + t^7}{x^2 + \sqrt 2\cdot x\cdot y\cdot t^4 + 2\cdot x\cdot t^8 - 2\cdot x + \frac{1}{2}\cdot y^2\cdot t^8 + \sqrt 2\cdot y\cdot t^{12} - \sqrt 2\cdot y\cdot t^4 + t^{16} - 2\cdot t^8 + 1}\\
      t^2
    \end{pmatrix}
    \endgroup
  \]
  As ultimate check we confirmed that $\kappa_{OY} = \kappa_{OY} \circ \iota_{OY}$. These computations are carried out 
  in the file \\
  \verb|enriques_projection.jl|.
\end{proof}

\subsection{Specialisation to positive characteristic}\label{sec:specialization}
In practice, it is computationally challenging and thus inconvenient to work over the field of rational numbers and even more so over a number field of degree $8$. Therefore we specialised the equations to positive characteristic when experimenting. The general theory is as follows:\\

Let $X_K$ be K3 surface over a number field $K$ and let $\mathfrak{p}$ be a finite place of $K$ where $X_K$ has good reduction. Choose a smooth model
of $X_K$ over the localized ring of integers of $\mathcal O_{K, \mathfrak{p}}$ and let $X_\mathfrak{p}$ be its special fiber.
Let $X$ denote the base change of $X_K$ to the algebraic closure $\overline{K}$ and $X_{\overline{k}}$ the base change of $X_{\mathfrak{p}}$ to the algebraic closure of the residue field $k$.
Specialisation of (effective) divisors induces a specialisation map
\begin{equation}\label{eqn:specialization}
\mathrm{sp}: \NS(X) \to \NS(X_{\overline{k}}),
\end{equation}
which is injective \cite[Proposition 6.2]{luijk} and preserves intersection products \cite[Corollary 20.3]{Fulton98}. Note that in general it fails to be surjective \cite{Charles14} \cite{CostaElsenhansJanel}.
Under the additional hypothesis that $\mathrm{sp}$ is surjective and $\mathrm{sp}(\mathcal{L})$ is ample for at least one ample line bundle $\mathcal{L}$ on $X$, the specialisation map $\mathrm{sp}$ is equivariant with respect to the group homomorphism $\Aut(X)\to \Aut(X_{\overline{k}})$ \cite[Theorem 2.1]{lieblich-maulik} and the latter is injective.\\

Let us return to our surface $X_{e_1}/\QQ$ constructed in the previous section.
The choice of the Weierstrass equation $y^2=x^3+(1-t^8)x$ defines a model of $X_{e_1}$ over $\Spec \mathbb{Z}$ which has good reduction away from $2$.
Recall that to access $\NS(X_{e_1})$, we had to perform a base change to the number field $K/\QQ$ of degree $8$.
In order to compute over a prime field, one has to choose an odd prime $\mathfrak{p} \leq \mathcal{O}_K$ of degree $1$. We picked a prime over $113$ ($41$ works just as well), so that the residue field $k$ equals $\FF_{113}$.
Finally, note that if $f$ is an automorphism of $X_{e_1}/\mathcal{O}_K$ and $\overline{f} \in \Aut(X_k)$ its specialisation, then $\mathrm{sp} \circ f^* = \overline{f}^* \circ \mathrm{sp}$.
These considerations pay off in practice by the fact that both, the correctness of the final result, as well as 
some intermediate computations of intersection numbers, can be certified in positive characteristic.

For our practical computations, however, it is not necessary to define a model 
of $E_1$ over $\ZZ$, a ring of algebraic integers, or a localization thereof 
as a data structure of its own. Instead, we will be working with $E_1$ over 
the number field $K$ introduced above and utilize the projection map
\begin{equation}
  \label{eqn:GoodReductionMapFields}
  K \dashrightarrow \FF_{113}, \quad \frac{p}{q} \mapsto (p + \mathfrak{p})\cdot (q + \mathfrak{p})^{-1} \quad 
  \textnormal{for} \quad q \notin \mathfrak{p}
\end{equation}
which is defined only on the localization $\OO_{K, \mathfrak{p}} \subset K$.
Should this map ever be applied to a non-admissible element, an error will 
be thrown to indicate that the curve or scheme's concrete model is not defined over 
the ring $\OO_{K, \mathfrak{p}}$.

From an experimental point of view, it is much more convenient to continue with
$\overline X_{e_1}$ and search for the specialization
$\overline f_{OY} \colon (X_{e_1})_k \to (X_{e_1})_k$ of the
automorphism $f_{OY}$ in characteristic zero. This is what we 
did in practice. After we found an automorphism $\overline f_{OY}$ 
with the anticipated action on 
$\mathrm{NS}((X_{e_1})_k)$, we let it guide our search for $f_{OY}$ and
matched the objects in characteristic zero to their reductions over $k=\FF_{113}$.
Only this final calculation is visible and carried out in 
\verb|minimum_entropy_automorphism.jl|.

\subsection{Fixing generators of the N\'eron-Severi lattice}
The \textit{trivial lattice} of $X_{e_1}$ is automatically extracted in
\Oscar from the resolution of singularities of the Weierstrass
model.  In  general, the computation of the Mordell-Weil group of an elliptic curve  is difficult and there is no automatism in \Oscar. However, by Remark \ref{rem:FurtherComputationsInFirstConstruction}, we have found generators of a sublattice $L$ of 
index $2$ in $\NS(X_{e_1})$, consisting of fiber components and sections. These can be transferred to the Weierstrass model. We did the pushforward in positive characteristic.
To also find the sections in characteristic zero, i.e. the elements of the Mordell-Weil
group, we apply a multivariate Newton iteration to $\FF_{113}(t)$-rational points on the reduction $\overline E_1$ over $\FF_{113}(t)$ to obtain $\QQ_{113}(t)$-rational points of $E_1$ (up to finite precision in the $113$-adic numbers $\QQ_{113}$). These are then used to reconstruct $K(t)$-rational points on $E_1$ via lattice reduction.
Note that since the map $\mathrm{sp}$ in Equation \ref{eqn:specialization} is injective, and the curves considered are unique in their numerical equivalence class, the lift is unique.

The Mordell-Weil
group of $\overline E_1$ is generated by the points with coordinates $(x, y)$ given by
\begin{equation}
  \label{eqn:RationalPointsPosChar}
  \begin{matrix}
    \overline{S_1} &=&  (62 t^4 + 6 t^3 + 106 t^2 + 75 t + 26,& 66 t^6 + 67 t^5 + 37 t^4 + 60 t^3 + 67 t^2 + 101 t + 58)\\
 \overline{S_2} &=&(112 t^4 + 97 t^2 + 98,& 104 t^5 + 82 t^3 + 91 t)\\
 \overline{S_3} &=&(t^4 + 112,& 87 t^4 + 26)\\
 \overline{S_4} &=&(8 t^4 + 17 t^3 + 11 t^2 + 70 t + 7,& 39 t^6 + 17 t^5 + 60 t^4 + 64 t^3 + 41 t^2 + 29 t + 89)\\
 \overline{S_5} &=&(50 t^4 + 104 t^2 + 72,& 56 t^6 + 30 t^4 + 36 t^2 + 104)\\
 \overline{S_6} &=&(31 t^4 + 45 t^3 + 60 t^2 + 59 t + 13,& 43 t^6 + 18 t^5 + 67 t^4 + 75 t^3 + 10 t^2 + 44 t + 17)\\
 \overline{S_7} &=&(105 t^4 + 43 t^3 + 61 t^2 + 17 t + 7,& 20 t^6 + 33 t^5 + 53 t^4 + 9 t^3 + 63 t^2 + 70 t + 89)\\
 \overline{S_8} &=&(62 t^4 + 75 t^3 + 8 t^2 + 107 t + 87,& 66 t^6 + 103 t^5 + 103 t^4 + 63 t^3 + 46 t^2 + 37 t + 34)\\
  \end{matrix}
\end{equation}
Liftings $S_i$ of these points to $E_1/K(t)$ are given in Figure \ref{fig:LiftingsToCharZero}.
We let $s_i = [S_i] \in \NS(X_{e_1})$ for $i=1,\dots, 8$.\\

Let $o_1$ be the zero section, and $e_1$ the class of a fiber. By Tate's algorithm \cite[IV \S9 Table 4.1]{Silverman94}, we have eight reducible fibers at $t^8=1$ consisting of two $(-2)$ curves meeting tangentially in one point. These are said to be of Kodaira type III and form and $\widetilde{A}_1$-configuration. Each one of these fibers provides us with a single component $a_{1,i}$ which does
not meet $o_1$. We fix the indices $i=1,\dots, 8$ in such a way that $a_{1,1},\dots a_{8,1}$ are contained in fibers over the points
\[\arraycolsep=2pt\def\arraystretch{2}
\begin{array}{lccclcc}
  t_1 &=& -1, & \quad &
t_2 &=&-\frac{1}{4}\alpha^6 - \frac{1}{4}\alpha^4 - \frac{7}{4}\alpha^2 - \frac{3}{4}, \\
  t_3 &=&\frac{1}{4}\alpha^6 - \frac{1}{4}\alpha^4 + \frac{7}{4}\alpha^2 - \frac{3}{4}, & \quad &
t_4 &=&-\frac{1}{2}\alpha^6 - \frac{5}{2}\alpha^2, \\
  t_5 &=&\frac{1}{2}\alpha^6 + \frac{5}{2}\alpha^2, & \quad &
t_6 &=&-\frac{1}{4}\alpha^6 + \frac{1}{4}\alpha^4 - \frac{7}{4}\alpha^2 + \frac{3}{4}, \\
  t_7 &=&\frac{1}{4}\alpha^6 + \frac{1}{4}\alpha^4 + \frac{7}{4}\alpha^2 + \frac{3}{4}, & \quad &
t_8 &=& 1.
\end{array}
\]
where $\alpha$ is fixed with $\alpha^8+6\alpha^4+1=0$.

Let $B = (e_1,o_1,a_{1,1},\dots a_{1,8},s_1,\dots s_8)$. The gram matrix with respect to the basis $B$ of $\NS(X_{e_1})\otimes \QQ$ is
\[G_B=\left(\begin{array}{rrrrrrrrrrrrrrrrrr}
0 & 1 & 0 & 0 & 0 & 0 & 0 & 0 & 0 & 0 & 1 & 1 & 1 & 1 & 1 & 1 & 1 & 1 \\
1 & -2 & 0 & 0 & 0 & 0 & 0 & 0 & 0 & 0 & 0 & 0 & 0 & 0 & 0 & 0 & 0 & 0 \\
0 & 0 & -2 & 0 & 0 & 0 & 0 & 0 & 0 & 0 & 1 & 0 & 1 & 0 & 1 & 1 & 1 & 0 \\
0 & 0 & 0 & -2 & 0 & 0 & 0 & 0 & 0 & 0 & 0 & 1 & 0 & 0 & 0 & 0 & 1 & 1 \\
0 & 0 & 0 & 0 & -2 & 0 & 0 & 0 & 0 & 0 & 1 & 0 & 0 & 1 & 1 & 0 & 1 & 1 \\
0 & 0 & 0 & 0 & 0 & -2 & 0 & 0 & 0 & 0 & 1 & 1 & 1 & 0 & 0 & 0 & 1 & 0 \\
0 & 0 & 0 & 0 & 0 & 0 & -2 & 0 & 0 & 0 & 0 & 1 & 1 & 1 & 0 & 1 & 0 & 1 \\
0 & 0 & 0 & 0 & 0 & 0 & 0 & -2 & 0 & 0 & 0 & 0 & 0 & 0 & 1 & 1 & 0 & 0 \\
0 & 0 & 0 & 0 & 0 & 0 & 0 & 0 & -2 & 0 & 1 & 1 & 0 & 1 & 0 & 1 & 0 & 0 \\
0 & 0 & 0 & 0 & 0 & 0 & 0 & 0 & 0 & -2 & 0 & 0 & 1 & 1 & 1 & 0 & 0 & 1 \\
1 & 0 & 1 & 0 & 1 & 1 & 0 & 0 & 1 & 0 & -2 & 2 & 2 & 2 & 2 & 0 & 0 & 2 \\
1 & 0 & 0 & 1 & 0 & 1 & 1 & 0 & 1 & 0 & 2 & -2 & 1 & 0 & 2 & 2 & 1 & 1 \\
1 & 0 & 1 & 0 & 0 & 1 & 1 & 0 & 0 & 1 & 2 & 1 & -2 & 0 & 0 & 1 & 2 & 1 \\
1 & 0 & 0 & 0 & 1 & 0 & 1 & 0 & 1 & 1 & 2 & 0 & 0 & -2 & 1 & 2 & 2 & 1 \\
1 & 0 & 1 & 0 & 1 & 0 & 0 & 1 & 0 & 1 & 2 & 2 & 0 & 1 & -2 & 1 & 1 & 1 \\
1 & 0 & 1 & 0 & 0 & 0 & 1 & 1 & 1 & 0 & 0 & 2 & 1 & 2 & 1 & -2 & 2 & 2 \\
1 & 0 & 1 & 1 & 1 & 1 & 0 & 0 & 0 & 0 & 0 & 1 & 2 & 2 & 1 & 2 & -2 & 1 \\
1 & 0 & 0 & 1 & 1 & 0 & 1 & 0 & 0 & 1 & 2 & 1 & 1 & 1 & 1 & 2 & 1 & -2
\end{array}\right).
\]
The lattice $\NS(X_{e_1})$ is generated by $B$ together with the $2$-torsion section $p_2$.
The latter has the following coordinates with respect to $B$.
\[
v_B(p_2)=\frac{1}{2}\left(\begin{array}{rrrrrrrrrrrrrrrrrr}
4 & 2 & -1 & -1 & -1 & -1 & -1 & -1 & -1 & -1 & 0 & 0 & 0 & 0 & 0 & 0 & 0 & 0
\end{array}\right).
\]
The $G$-invariant ample class is given by
\[v_B(h)=
\frac{1}{2}\left(\begin{array}{rrrrrrrrrrrrrrrrrr}8& 2& -1& -1& -1& 1& -1& -1& 1& -1& 4& 4& 0& -2& 2& -2& -4& 0
\end{array}\right).\]
Let $\epsilon$ denote the involution on $X_{e_1}$ given in Weierstrass coodinates by $(x,y,t) \mapsto (x,y,-t)$.
By direct computations we find that the Enriques involution $j \in G$ of $X_h$ is given on $X_{e_1}$ by $T_P \circ \epsilon$ where
$T_P$ denotes the translation automorphism induced by the section $P$ with coordinates
\begin{eqnarray*}
  x &=& \frac{(-\alpha^6 - \alpha^4 - 7\alpha^2 + 1)t + 4}{4t + \alpha^6 - \alpha^4 + 7\alpha^2 + 1},\\
  y &=& \frac{2}{2t^3 + (\alpha^6 - \alpha^4 + 7\alpha^2 + 1)t^2 + (-\alpha^4 + \alpha^2 - 1)t} \cdot \left(
  \left(\frac{3}{2}\alpha^7 - \alpha^5 + \frac{19}{2}\alpha^3 - 8\alpha\right)t^4\right. \\
  & & + \left(-10\alpha^7 - 2\alpha^5 - 58\alpha^3 - 18\alpha\right)t^3 + \left(-12\alpha^7 + 6\alpha^5 - 72\alpha^3 + 30\alpha\right)t^2 \\
  & & \left.+ \left(8\alpha^7 + 4\alpha^5 + 44\alpha^3 + 24\alpha\right)t + \frac{7}{2}\alpha^7 - \alpha^5 + \frac{39}{2}\alpha^3 - 4\alpha \right).
\end{eqnarray*}
To obtain a simpler Enriques involution on $X_{e_1}$ compose $\epsilon$ with the translation by the unique $2$-torsion section. The resulting involution $\iota$ is an Enriques involution and $X/\iota$ is an Enriques surface of base-change type. Using our \Oscar implementation, we compute the coinvariant lattice $\NS(X_{e_1})_\iota$ and check that it is in the correct genus. Then Proposition \ref{prop:IdentifyInvolution} yields that $\iota$ is $\Aut(Y)$-conjugate to $\iota_{OY}$ and that $X_{e_1}/\iota:=Y_{f_1} \cong Y_{OY}$.\\

It remains to express $(f_{OY})_*$ and $(\iota_{OY})_*$ in terms of the basis $B$ of $\NS(X_{e_1})$. So far we have it only as an isometry on the lattice $L_{K3}$ coming from the database entry \verb|40.8.0.1|.
It comes with an implicitly given Hodge structure and $G$-invariant ample class $h$.

We thus have to define an embedding $\gamma \colon \NS(X_{h}) \to L_{K3}$ that can be extended to a Hodge isometry $H^2(X_{h},\ZZ)\to L_{K3}$.
Recall from Proposition \ref{prop:rational_curves_degrees} that in $L_{K3}$ we are given $64 = 32+32$ vectors $(c_1,\dots, c_{64})=C$ of square $(-2)$. They correspond to $64$ explicit classes $D=(d_1,\dots, d_{64})$ spanning an index $2$ sublattice of $\NS(X_{h})$. One can view $C$ and $D$ as graphs with edge-weights, compute a graph isomorphism\footnote{We used sageMath \cite{sagemath} for this purpose.} 
$\delta \colon C \to D$
and define $\gamma$ as its linear extension (if it exists). Then, one has to modify the isomorphism in such a way that it extends to a Hodge isometry.
Finally, one computes the pushforward $\NS(X_{h}) \to \NS(X_{e_1})$, to represent everything with respect to the basis $B$ of the elliptic fibration on $X_{e_1}$.

A more direct way to go forward is the following:
We compute the image of $\Aut(Y_{f_1}) \to O(\Num(Y_{f_1}))$ and use it to find an automorphism of the correct entropy.
Since the intermediate results and details of these computations are not of relevance, we just present the end result and how it can be checked for correctness.

The matrix representation of the pushforward $(f_{OY})_*$ on $\NS(X_{e_1})$ with respect to the basis $B$ is
\begin{equation}\label{eqn:fOY_pushforward}
D_{B}((f_{OY})_*) = \left(\begin{array}{cccccccccccccccccc}
3 & 1 & -1 & -1 & 1 & -1 & 0 & -1 & 1 & 0 & 1 & 0 & -1 & 1 & 0 & 0 & 0 & 0 \\
2 & 1 & 0 & 0 & 0 & 0 & -1 & -1 & 0 & 0 & 1 & 0 & 0 & 0 & 0 & -1 & 0 & 0 \\
3 & 1 & -1 & -1 & 1 & -1 & -1 & -1 & 1 & 0 & 1 & 0 & -1 & 1 & 0 & 0 & 0 & 0 \\
2 & 1 & -1 & 0 & 1 & -1 & 0 & -1 & 1 & 0 & 1 & 0 & -1 & 1 & 0 & 0 & 0 & 0 \\
3 & 1 & -1 & -1 & 0 & -1 & 0 & -1 & 1 & 0 & 1 & 0 & -1 & 1 & 0 & 0 & 0 & 0 \\
2 & 1 & 0 & -1 & 1 & -1 & 0 & -1 & 1 & 0 & 1 & 0 & -1 & 1 & 0 & 0 & 0 & 0 \\
0 & 0 & 0 & 0 & 0 & 0 & 0 & 0 & 0 & 1 & 0 & 0 & 0 & 0 & 0 & 0 & 0 & 0 \\
1 & 0 & 0 & 0 & 0 & 0 & 0 & -1 & 0 & 0 & 0 & 0 & 0 & 0 & 0 & 0 & 0 & 0 \\
0 & 0 & 0 & 0 & 0 & 0 & 0 & 0 & 1 & 0 & 0 & 0 & 0 & 0 & 0 & 0 & 0 & 0 \\
1 & 0 & 0 & 0 & 0 & -1 & 0 & 0 & 0 & 0 & 0 & 0 & 0 & 0 & 0 & 0 & 0 & 0 \\
-3 & -2 & \frac{1}{2} & \frac{1}{2} & \frac{3}{2} & \frac{1}{2} & \frac{1}{2} & -\frac{1}{2} & \frac{1}{2} & \frac{1}{2} & 1 & 0 & 0 & 1 & 0 & 0 & 1 & 0 \\
14 & 6 & -\frac{5}{2} & -\frac{5}{2} & -\frac{1}{2} & -\frac{3}{2} & -\frac{5}{2} & -\frac{5}{2} & -\frac{1}{2} & -\frac{3}{2} & 2 & 0 & -1 & 0 & 0 & -2 & -1 & -1 \\
8 & 4 & -2 & -1 & 0 & -1 & -1 & -2 & 0 & -1 & 1 & 0 & -1 & 0 & 0 & -1 & 0 & 0 \\
6 & 3 & -\frac{3}{2} & -\frac{3}{2} & \frac{3}{2} & -\frac{3}{2} & -\frac{3}{2} & -\frac{3}{2} & \frac{1}{2} & -\frac{1}{2} & 2 & 0 & -2 & 1 & 1 & -1 & -1 & 0 \\
10 & 5 & -2 & -1 & 0 & -1 & -1 & -2 & 0 & -1 & 1 & 0 & -1 & 0 & -1 & -1 & 0 & 0 \\
2 & 1 & -\frac{1}{2} & -\frac{1}{2} & \frac{1}{2} & -\frac{1}{2} & \frac{1}{2} & -\frac{3}{2} & \frac{1}{2} & -\frac{1}{2} & 1 & 0 & 0 & 1 & -1 & 0 & 1 & 0 \\
0 & -1 & 0 & 0 & 1 & 1 & 0 & -1 & 1 & 0 & 2 & 1 & 0 & 0 & 0 & -1 & 0 & 0 \\
10 & 4 & -\frac{5}{2} & -\frac{5}{2} & \frac{5}{2} & -\frac{3}{2} & -\frac{3}{2} & -\frac{7}{2} & \frac{3}{2} & -\frac{1}{2} & 3 & 0 & -2 & 2 & 0 & -2 & 0 & 0
\end{array}\right)
\end{equation}
and the pushforward of the accompanying Enriques involution reads
\begin{equation}\label{eqn:iOY_pushforward}
 D_B(\iota_*) = \left(
 \begin{array}{cccccccccccccccccc}
1 & 0 & 0 & 0 & 0 & 0 & 0 & 0 & 0 & 0 & 0 & 0 & 0 & 0 & 0 & 0 & 0 & 0 \\
2 & 1 & -\frac{1}{2} & -\frac{1}{2} & -\frac{1}{2} & -\frac{1}{2} & -\frac{1}{2} & -\frac{1}{2} & -\frac{1}{2} & -\frac{1}{2} & 0 & 0 & 0 & 0 & 0 & 0 & 0 & 0 \\
1 & 0 & 0 & 0 & 0 & 0 & 0 & 0 & 0 & -1 & 0 & 0 & 0 & 0 & 0 & 0 & 0 & 0 \\
1 & 0 & 0 & 0 & 0 & 0 & 0 & 0 & -1 & 0 & 0 & 0 & 0 & 0 & 0 & 0 & 0 & 0 \\
1 & 0 & 0 & 0 & 0 & 0 & 0 & -1 & 0 & 0 & 0 & 0 & 0 & 0 & 0 & 0 & 0 & 0 \\
1 & 0 & 0 & 0 & 0 & 0 & -1 & 0 & 0 & 0 & 0 & 0 & 0 & 0 & 0 & 0 & 0 & 0 \\
1 & 0 & 0 & 0 & 0 & -1 & 0 & 0 & 0 & 0 & 0 & 0 & 0 & 0 & 0 & 0 & 0 & 0 \\
1 & 0 & 0 & 0 & -1 & 0 & 0 & 0 & 0 & 0 & 0 & 0 & 0 & 0 & 0 & 0 & 0 & 0 \\
1 & 0 & 0 & -1 & 0 & 0 & 0 & 0 & 0 & 0 & 0 & 0 & 0 & 0 & 0 & 0 & 0 & 0 \\
1 & 0 & -1 & 0 & 0 & 0 & 0 & 0 & 0 & 0 & 0 & 0 & 0 & 0 & 0 & 0 & 0 & 0 \\
2 & 1 & -1 & 0 & 0 & -1 & 0 & 0 & 0 & 0 & 0 & 0 & -1 & 1 & 0 & 0 & 0 & 0 \\
4 & 2 & -\frac{1}{2} & -\frac{1}{2} & -\frac{1}{2} & -\frac{1}{2} & -\frac{1}{2} & -\frac{1}{2} & -\frac{1}{2} & -\frac{1}{2} & 0 & -1 & 0 & 0 & 0 & 0 & 0 & 0 \\
0 & 0 & \frac{1}{2} & -\frac{1}{2} & -\frac{1}{2} & \frac{1}{2} & \frac{1}{2} & -\frac{1}{2} & -\frac{1}{2} & \frac{1}{2} & 0 & 0 & 1 & 0 & 0 & 0 & 0 & 0 \\
-2 & -1 & 1 & 0 & 0 & 1 & 0 & 0 & 0 & 0 & 1 & 0 & 1 & 0 & 0 & 0 & 0 & 0 \\
0 & 0 & \frac{1}{2} & -\frac{1}{2} & \frac{1}{2} & -\frac{1}{2} & -\frac{1}{2} & \frac{1}{2} & -\frac{1}{2} & \frac{1}{2} & 0 & 0 & 0 & 0 & 1 & 0 & 0 & 0 \\
2 & 1 & -1 & 0 & 0 & -1 & 0 & 0 & 0 & 0 & -1 & 0 & -1 & 1 & 0 & 1 & 0 & 0 \\
4 & 2 & -1 & -1 & 0 & -1 & -1 & 0 & 0 & 0 & 1 & 0 & -1 & 0 & 1 & -1 & -1 & 0 \\
12 & 6 & -1 & -1 & -3 & 0 & -2 & -1 & -2 & -2 & -1 & 0 & 1 & -2 & -1 & -1 & 0 & -1
\end{array}
 \right).
\end{equation}
At this point one should perform a number of sanity checks on the matrix $D_{B}((f_{OY})_*)$:
\begin{enumerate}
\item it is an isometry, i.e. $D_{B}((f_{OY})_*) G_B D_{B}((f_{OY})_*)^T = G_B$
\item its spectral radius is $\tau_8$,
\item it preserves the nef cone, i.e. $h':=(f_{OY})_*(h)$ is ample,
\item it can be extended to a Hodge isometry,
\item it commutes with $D_{B}(\iota_*)$.
\end{enumerate}
(1) and (2) are trivial to check using any computer algebra system.
For (3) one certifies that the set
\[\{r \in \NS(X_{e_1}) \mid r^2=-2, r.h>0, r.h'<0\}\]
is empty. It is computed using the \Oscar function \verb|separating_hyperplanes| based on \cite[\S 3]{Shimada24}.
For (4) one computes its action on the discriminant group and checks that it generates the same subgroup of $O(\NS(X)^\vee/ \NS(X))$, as the action of the non-symplectic automorphism of order $8$ given by $t\mapsto \zeta_8 t$.
These checks are carried out in the file \verb|sanity_check.jl|.

\section{Construction of the automorphism}\label{sec:ConstructionAutomorphism}

In order to find $f_{OY}$, we apply techniques called
``Kneser's neighbor method'' and ``fibration hopping''. 
These have been described in 
\cite{Kumar14} and \cite{ElkiesKumar14} and they have already been 
put to use in 
\cite{ElkiesSchuett15}, \cite{Kumar15}, and \cite{BrandhorstElkies23}.
As it turned out, we do not need to apply these methods in full generality 
to solve our particular problem and we shall henceforth confine ourselves
to describing the minimal requirements only. 

Given that we already know the prescribed action of $f_{OY}$ on the 
N\'eron-Severi lattice $\mathrm{NS}(X_{e_1})$ of our K3 surface, 
Equation \ref{eqn:fOY_pushforward},
 we may consider the image
of the class of a fiber of $\pi_1$
\[
  e_1' := (f_{OY})_*(e_1) \in \mathrm{NS}(X_{e_1}).
\]
Its basis representation is given by the first row of the matrix in Equation \ref{eqn:fOY_pushforward}.
Since $f_{OY}$ is an automorphism of the very same surface, $e'$ appears as the class of a fiber of another elliptic fibration
\[
  \pi' \colon X_{e_1'} \to \PP^1
\]
which must be isomorphic to the original one by means of $f_{OY}$. Our next goal is to explicitly compute the elliptic fibration $\pi'$. This can be done by fibration hopping, which in a nutshell is the following.\\

Recall that elliptic fibrations on a K3 surface $X$ are in bijection with primitive, isotropic and nef classes $e \in \NS(X)$. Such classes are called \textit{elliptic}. Indeed the linear system $|e|$ of an elliptic class induces an elliptic fibration (possibly without section) and conversely, the class of a fiber of any elliptic fibration is an elliptic class. Two elliptic fibrations corresponding to elliptic classes $e,e' \in \NS(X)$ are called $n$-\textit{neighbors} if $e.e'=n$.

Given an elliptic fibration $\phi \colon X \to \PP^1$ with elliptic class $e$ and
an elliptic class $e' \in \mathrm{NS}(X)$ with $e.e'=2$, i.e. the class of a $2$-neighborly fibration, a \textit{two neighbor step} is a procedure to compute the fibration $\phi'\colon X \to \mathbb{P}^1$ corresponding to $e'$ in terms of a weierstrass model and the corresponding changes of coordinates.

In our case $e_1.e_1'=2$ and so $\phi:=\pi_1$ and $\phi':=\pi'$ are $2$-neighbors.
\subsection{A $2$-neighbor step}
We proceed as follows. First we compute the linear system of $e'$ by means of the methods
outlined in \cite[Lemma 2.14]{BrandhorstElkies23}. This provides us with an elliptic-fibration $\pi' \colon X \to \PP^1$, described by a rational function
on $X_{e_1}$, the pullback of an affine coordinate of $\PP^1$.

A well established procedure (which has, for example, also been described and 
used in \cite{BrandhorstElkies23}) allows to extract a new equation of the form
\[
  \tilde y^2 = g(\tilde x), \quad g(\tilde x) \in K(\tilde t)[\tilde x] \textnormal{ of degree } \leq 4,
\]
from the original Weierstrass equation (\ref{eqn:ModifiedFirstWeierstrassEquation}), 
together with the rational change of coordinates
\[
  \Phi \colon K(\tilde t, \tilde x, \tilde y) \to K(t, x, y)
\]
taking (\ref{eqn:ModifiedFirstWeierstrassEquation}) to $\tilde y^2 - g(\tilde x)$.
This has been automatized in the function \verb|two_neighbor_step| 
called in \newline\verb|minimum_entropy_automorphism.jl|.

Next, we need to find a section for $\pi'$, i.e. a $K(\tilde t)$-rational point.
Finding rational points is known to be a hard problem. Hence, it still
has to be carried out by hand and constitutes the main obstruction to
full automatization of $2$-neighbor steps.

We are lucky since the section $[D]:=s_8$ of $\pi$ satisfies $s_8.e'=1$. Hence its pushforward $\Phi_*D$ induces a $K(\tilde t)$-rational point.
However, computations for the extraction of explicit coordinates as a point 
on the curve defined by $\tilde y^2 - g(\tilde x)$ do not terminate in \Oscar at the moment.
For the reduction in positive characteristic computations do terminate in \Oscar and 
we can match the result with the reduction of the hardcoded point in characteristic zero. 

This information, in turn, provides the input to a constructor for an elliptic surface
which automatically takes care of transforming the given equation into a globally minimal integral Weierstrass
form and returns the elliptic surface, together with the additional rational
transformation. The final Weierstrass equation defining our second model 
$X_{e'}$ for this elliptic surface is
\begin{eqnarray*}
  \tilde y^2 &=& \tilde x^3 - ((963\alpha^7 + 468\alpha^6 + 27\alpha^5 - 168\alpha^4 + 5613\alpha^3 + 2724\alpha^2 + 165\alpha - 987)\tilde t^8 \\
  & & + (8146\alpha^7 + 8200\alpha^6 + 4034\alpha^5 + 328\alpha^4 + 47486\alpha^3 + 47720\alpha^2 + 23654\alpha + 1776)\tilde t^7 \\
  & & + (-8946\alpha^7 + 29456\alpha^6 + 30142\alpha^5 + 15750\alpha^4 - 52094\alpha^3 + 171080\alpha^2 + 176890\alpha + 90594)\tilde t^6 \\
  & & + (-198380\alpha^7 - 23604\alpha^6 + 59304\alpha^5 + 66080\alpha^4 - 1156092\alpha^3 - 140588\alpha^2 + 351904\alpha + 378672)\tilde t^5 \\
  & & + (-505540\alpha^7 - 269710\alpha^6 - 41160\alpha^5 + 81760\alpha^4 - 2946020\alpha^3 - 1581790\alpha^2 - 219520\alpha + 455420)\tilde t^4 \\
  & & + (-376040\alpha^7 - 426636\alpha^6 - 243656\alpha^5 - 29148\alpha^4 - 2190552\alpha^3 - 2506924\alpha^2 - 1378328\alpha - 212716)\tilde t^3 \\
  & & + (129024\alpha^7 - 187768\alpha^6 - 238056\alpha^5 - 118762\alpha^4 + 753592\alpha^3 - 1120308\alpha^2 - 1334368\alpha - 746354)\tilde t^2 \\
  & & + (232488\alpha^7 + 47498\alpha^6 - 62932\alpha^5 - 62418\alpha^4 + 1356168\alpha^3 + 258014\alpha^2 - 328156\alpha - 403254)\tilde t \\
  & & + 53663\alpha^7 + 32688\alpha^6 + 4351\alpha^5 - 4928\alpha^4 + 313121\alpha^3 + 184520\alpha^2 + 37697\alpha - 41349)\cdot \tilde x
\end{eqnarray*}
The rational transformations produced by this process can be used to construct an explicit
identification 
\[
  \Xi \colon X_{e'} \overset{\cong}{\longrightarrow} X_{e_1}
\]
of our two instances of 
the K3 surface arising from the different Weierstrass models. This is what we call 
a \textit{geometric realization} of a $2$-neighbor step.

\begin{remark}
  \label{rem:NonUniquenessOfProcedure}
  Note that the identification $\Xi$ is not unique! It can, for instance be
  modified by the \textit{auto}morphism of $X_{e'}$ given by translation with
  a section, or the one induced by an admissible  M\"obius transformation
  (see Definition \ref{def:AdmissibleMoebiusTransformation}) on
  the base $\PP^1$ of the elliptic fibration. In fact, this ambiguity
  will become important below.

  Moreover, these ambiguities and the related automated choices introduce 
  a certain randomness in the intermediate outputs of our construction script
  \verb|minimum_entropy_automorphism.jl|.
  It is important to note that therefore 
  many of the following matrices and other intermediate results 
  may differ in a reproduction of our findings; we nevertheless leave them 
  here for illustration and comparison. 
\end{remark}

\subsection{Completion to an automorphism of $X_{e_1}$}
Up to this point we produced the identification $\Xi \colon X_{e'} \to X_{e_1}$.
Since the elliptic fibrations $\pi \colon X_{e_1} \to \PP^1$ and
$\pi' \colon X_{e'} \to \PP^1$ are \textit{abstractly isomorphic}, we can,
on the one hand, think 
of $h$ as a ``change of (Weierstrass-) coordinates'' on the same surface. 
On the other hand, we may \textit{choose} an identification $\Psi$ 
of these two different descriptions of this fibration. 
\begin{equation}
  \label{eqn:FibrationHoppingDiagram1}
  \begin{xy}
    \xymatrix{
      X_{e_1} \ar[d]_{\pi} &
      X_{e'} \ar[d]_{\pi'} \ar[l]^\cong_{\Xi}&
      X_{e_1} \ar[l]_{\Psi} \ar[d]_{\pi} \\
      \PP^1 &
      \PP^1 &
      \PP^1 \ar[l]_{\mu}
    }
  \end{xy}
\end{equation}
The composition $\Xi \circ \Psi \colon X_{e_1} \to X_{e_1}$ is then already an honest automorphism
of $X_{e_1}$; the latter also considered as a single instance of this object in memory.

\medskip
To get our hands on possible choices for $\Psi$, we implemented the internal function \newline
\verb|admissible_moebius_transformations| which takes two 
\verb|EllipticSurface|s as input. 

\begin{definition}
  \label{def:AdmissibleMoebiusTransformation}
  Given two elliptic fibrations $\pi, \pi' \colon X \to \PP^1$, 
  a M\"obius transformation $\mu \colon \PP^1 \to \PP^1$ is \textit{admissible}
  if it lifts to an automorphism $\Psi \colon X \to X$ such that 
  $\pi' \circ \Psi = \mu \circ \pi$.
\end{definition}

Consequently, the function computes a full set of those M\"obius 
transformations of the base $\PP^1$ which take the configuration of critical 
values of the first elliptic fibration into the configuration of the second. Note that 
this list might already be empty, but it is certainly finite once we have at least three 
such critical values. 

In a second step, an attempt is made to lift every such transformation 
$\mu \colon \PP^1 \to \PP^1$ to an isomorphism of the elliptic surfaces. 
Whether or not such an identification exists is decided from the 
generic fibers after base change by $\mu$ and in the affirmative case 
such an isomorphism is automatically extended to the associated elliptic 
surfaces.

\medskip
For such an arbitrary choice of $\Psi$ we then compute the pushforward of $\Xi \circ \Psi$ 
on $\mathrm{NS}(X_{e_1})$. This is done in a highly automated and optimized way via
the method \verb|pushforward_on_algebraic_lattices|, but uses parts of the 
generic machinery developed for algebraic cycles; see the next section on 
implementations for more details.

The explicit
form of $(\Xi \circ \Psi)_*$ then gives us a first impression on how far we still
are from $f_{OY}$. Note that, by construction, $\Xi \circ \Psi$ is an
automorphism of $X_{e_1}$ which takes the class of the fiber $F$ to the class of
$F' = (f_{OY})_*(F)$, but $\Xi \circ \Psi$ and $f_{OY}$ will in general be different
at this point.

\subsection{Resolving the final ambiguity}

To get from $\Xi\circ \Psi$ to $f_{OY}$ still requires some massaging. First
we may browse through the lifts of admissible M\"obius automorphisms for $X_{e_1}$ and
select one, $\Phi$ over $\nu$ say, for which $\Xi \circ \Psi \circ \Phi$ takes 
the reducible fibers of $X_{e_1}$ to the same ones as $f_{OY}$.
\begin{equation}
  \label{eqn:FibrationHoppingDiagram}
  \begin{xy}
    \xymatrix{
      X_{e_1} \ar[d]_{\pi} &
      X_{e'} \ar[d]_{\pi'} \ar[l]^\Xi_{\cong}&
      X_{e_1} \ar[l]_{\Psi} \ar[d]_{\pi}&
      X_{e_1} \ar[l]_\Phi \ar[d]_{\pi} \\
      \PP^1 &
      \PP^1 &
      \PP^1 \ar[l]_{\mu} &
      \PP^1 \ar[l]_{\nu}
    }
  \end{xy}
\end{equation}
Whether or not
this is the case can be read off from the matrices for the respective pushforwards 
on $\mathrm{NS}(X_{e_1})$. To this end we may inspect the pushforward of the composite
map $(f_{OY}^{-1})_* \circ (\Xi \circ \Psi)_* \circ \Phi_*$ on $\mathrm{NS}(X_{e_1})$,
which we can build from the already known information on the Hodge and lattice theoretic side.

Once we have found and chosen $\Phi$, the composite map $c = f_{OY}^{-1} \circ \Xi \circ \Psi \circ \Phi$
is an honest automorphism of the elliptic fibration $\pi \colon X_{e_1} \to \PP^1$, i.e. $\pi = \pi \circ c$.
At this point the representing matrix for the pushforward of this composite 
map $c$ on $\mathrm{NS}(X_{e_1})$ for the set of generators chosen earlier reads\footnote{cf. Remark \ref{rem:NonUniquenessOfProcedure}}
\[
  \begin{pmatrix}
\begin{array}{c|c|cccccccc|cccccccc}
    1& 0& 0& 0& 0& 0& 0& 0& 0& 0& 0& 0& 0& 0& 0& 0& 0& 0\\ 
    \hline
    -8 & -4 & \frac{1}{2} & \frac{1}{2} & \frac{5}{2} & -\frac{1}{2} & \frac{3}{2} & \frac{1}{2} & \frac{3}{2} & \frac{3}{2} & 1 & 0 & -1 & 2 & 1 & 1 & 0 & 1 \\
    \hline
    1 & 0 & -1 & 0 & 0 & 0 & 0 & 0 & 0 & 0 & 0 & 0 & 0 & 0 & 0 & 0 & 0 & 0 \\
0 & 0 & 0 & 1 & 0 & 0 & 0 & 0 & 0 & 0 & 0 & 0 & 0 & 0 & 0 & 0 & 0 & 0 \\
0 & 0 & 0 & 0 & 1 & 0 & 0 & 0 & 0 & 0 & 0 & 0 & 0 & 0 & 0 & 0 & 0 & 0 \\
1 & 0 & 0 & 0 & 0 & -1 & 0 & 0 & 0 & 0 & 0 & 0 & 0 & 0 & 0 & 0 & 0 & 0 \\
0 & 0 & 0 & 0 & 0 & 0 & 1 & 0 & 0 & 0 & 0 & 0 & 0 & 0 & 0 & 0 & 0 & 0 \\
1 & 0 & 0 & 0 & 0 & 0 & 0 & -1 & 0 & 0 & 0 & 0 & 0 & 0 & 0 & 0 & 0 & 0 \\
1 & 0 & 0 & 0 & 0 & 0 & 0 & 0 & -1 & 0 & 0 & 0 & 0 & 0 & 0 & 0 & 0 & 0 \\
0 & 0 & 0 & 0 & 0 & 0 & 0 & 0 & 0 & 1 & 0 & 0 & 0 & 0 & 0 & 0 & 0 & 0 \\
    \hline
-6 & -3 & \frac{1}{2} & \frac{1}{2} & \frac{3}{2} & -\frac{1}{2} & \frac{3}{2} & \frac{1}{2} & \frac{3}{2} & \frac{3}{2} & 0 & 0 & -1 & 2 & 1 & 1 & 0 & 1 \\
-5 & -3 & \frac{1}{2} & -\frac{1}{2} & \frac{5}{2} & -\frac{1}{2} & \frac{1}{2} & \frac{1}{2} & \frac{3}{2} & \frac{3}{2} & 1 & -1 & -1 & 2 & 1 & 1 & 0 & 1 \\
-5 & -3 & \frac{1}{2} & \frac{1}{2} & \frac{5}{2} & -\frac{1}{2} & \frac{1}{2} & \frac{1}{2} & \frac{3}{2} & \frac{1}{2} & 1 & 0 & -2 & 2 & 1 & 1 & 0 & 1 \\
-5 & -3 & \frac{1}{2} & \frac{1}{2} & \frac{3}{2} & -\frac{1}{2} & \frac{1}{2} & \frac{1}{2} & \frac{3}{2} & \frac{1}{2} & 1 & 0 & -1 & 1 & 1 & 1 & 0 & 1 \\
-5 & -3 & \frac{1}{2} & \frac{1}{2} & \frac{3}{2} & -\frac{1}{2} & \frac{3}{2} & \frac{1}{2} & \frac{3}{2} & \frac{1}{2} & 1 & 0 & -1 & 2 & 0 & 1 & 0 & 1 \\
-6 & -3 & \frac{1}{2} & \frac{1}{2} & \frac{5}{2} & -\frac{1}{2} & \frac{1}{2} & \frac{1}{2} & \frac{3}{2} & \frac{3}{2} & 1 & 0 & -1 & 2 & 1 & 0 & 0 & 1 \\
-5 & -3 & \frac{1}{2} & -\frac{1}{2} & \frac{3}{2} & -\frac{1}{2} & \frac{3}{2} & \frac{1}{2} & \frac{3}{2} & \frac{3}{2} & 1 & 0 & -1 & 2 & 1 & 1 & -1 & 1 \\
-4 & -3 & \frac{1}{2} & -\frac{1}{2} & \frac{3}{2} & -\frac{1}{2} & \frac{1}{2} & \frac{1}{2} & \frac{3}{2} & \frac{1}{2} & 1 & 0 & -1 & 2 & 1 & 1 & 0 & 0
\end{array}
\end{pmatrix}.
\]
Recall that with our choices made, the first $10$ generators span the trivial lattice
$\mathrm{Triv}(X_{e_1})$ while the remaining $8$ ones are sections of the elliptic fibration.
Note that the first row of the matrix above is the image of $e_1$. We see that $c_*(e_1)=e_1$ as expected.
For the class of a singular fiber $a_{1,i} + a_{2,i}$, we see in rows $3-10$ that $c_*(a_{1,i})$ is either $a_{1,i}$ or $a_{2,i}= e_1 - a_{1,i}$, confirming that $c_*$ acts trivially on $\PP^1$ (it has $8$ fixed points on $\PP^1$ below the reducible singular fibers).

The next step to resolve in the ambiguity for the choice 
of $\Psi$ are translations by sections.
The quotient $\mathrm{NS}(X_{e_1})/\mathrm{Triv}(X_{e_1})$ is isomorphic to the Mordell-Weil group of the generic fiber
$E_1$ of $X_{e_1}$ (cf. \cite[Theorem 6.5]{SchuettShioda19}).
A closer inspection of the second row of the above matrix, i.e. the image 
of the zero section $o$ on $X_{e_1}$, therefore yields that $o$ is taken to the section $p:=[P]$ with $P$ defined by the linear combination
\[P:=P_2 + S_1 - S_3 + 2S_4+S_5+S_6+S_8\]
in the group law of the generic fiber. Here the two torsion section $P_2$ comes in because of the denominator $2$ in the second row of the matrix.
Therefore, we compose with the translation $T$ by $-P$ on $X_{e_1}$.

\medskip
The very last ambiguity in the choice of $\Psi$ which is left at this point 
amounts to an automorphism of the generic fiber $E_1$.
In this case, there are four of them which we can get access to via \verb|automorphism_list|.
Their associated automorphisms of $X_{e_1}$ can be computed automatically via
the internal function \newline \verb|isomorphism_from_generic_fibers| and the induced 
maps on $\mathrm{NS}(X_{e_1})$ as usual with \newline \verb|pushforward_on_algebraic_lattices|.
It is now a simple matter of trial and error among the four automorphisms to find the correct final compensation
$A$.
By construction and since the action of $\Aut(X_{OY})$ on the N\'eron-Severi lattice is faithful (an element in the kernel would act trivially on $\NS(X_{e_1})^\vee/\NS(X_{e_1})\cong T(X_{OY})^\vee/T(X_{OY})$ and thus on $T(X_{e_1})$), we have constructed the automorphism of Oguiso and Yu as
\begin{eqnarray}
  \label{eqn:FinalAutomorphism}
  \mathrm{id} = A \circ T \circ f_{OY}^{-1} \circ \Xi \circ \Psi \circ \Phi &
  \Leftrightarrow & 
  f_{OY} = \Xi \circ \Psi \circ \Phi \circ A \circ T.
\end{eqnarray}
Since $(f_{OY})_*$ and $j_*$ commute, so do $f_{OY}$ and the enriques involution $j$.
Thus $f_{OY}$ descends to the Enriques surface $Y_{OY}$.

The final result can be checked for correctness by computing the pushforward $(f_{OY})_*$ explicitly, calculating its characteristic polynomial and checking that it divisible by the minimal polynomial $1-x^2-2x^3-x^4+x^6$ of $\tau_8$.

\section{On the implementation}\label{sec:Implementation}

The search for the automorphism of Oguiso and Yu relies heavily on computational tools, 
the development of which have been pushed in \Oscar with an eye towards this particular 
application. 
It is therefore of central interest for this paper to also report on these implementations. 

In particular, we have built up a framework for algebraic schemes\footnote{See also \cite{BrandhorstZachFruehbis} where we report on work in progress on the framework.} throughout the
roughly last three years as part of the project B5 of the SFB-TRR 195, Project-ID 286237555, 
by the DFG.
This joint effort has been carried out in wide parts
by M. Zach, with significant contributions from S. Brandhorst and A. Fr\"uhbis-Kr\"uger
and in exchange with J. B\"ohm. The experience from a previous implementation of
resolution of singularities in Singular \cite{Sin} and the efforts for parallelization 
in algebraic geometry from the joint work of A. Fr\"ubis-Kr\"uger and J. B\"ohm
\cite{BoehmFruehbisKrueger18} have been influential for the development.

The following is an overview on the cornerstones of these recent developments 
in \Oscar which are relevant for this particular publication. The discussion 
is by no means exhaustive and we make no attempt to prove all details with full 
mathematical rigor here, as this 
would expand the paper beyond any form acceptable for publication in a mathematical 
journal. Instead, we believe that the source code, which is Open Source and available online, 
can stand for itself, as it can be checked for mathematical rigor and correctness of 
implementation. For this note we refer to \Oscar version \verb|1.3.0-DEV| on commit
\verb|b4d5f32cbbcaa3bf68fdc798a3c4bc90ace3d758| in the official development branch of
\Oscar. 

\subsection{Polynomial rings, quotients, and localizations}\label{sec:algebras}
  Mathematically, an affine scheme of finite type over a ring $\bbk$ is of the form $X = \Spec R$ for some
  finitely generated $\bbk$-algebra $R$. The latter can be a polynomial ring over $\bbk$, or 
  a quotient thereof. By virtue of Rabinowitsch's trick this can be extended to 
  localizations 
  \[
    R[f^{-1}] \cong R[t]/\langle 1-t\cdot f\rangle
  \]
  at the multiplicative set given by the powers of an element $f \in R$. This is important as the spectra of such 
  localizations form a basis for the Zariski topology of $\Spec R$. 
  Therefore, a first step to realizing (affine) algebraic geometry in \Oscar 
  is to create the data types and functionality for 
  these rings:
  \begin{itemize}
    \item \verb|MPolyRing|: multivariate polynomial rings $P$ over a ring $\bbk$;
    \item \verb|MPolyQuoRing|: quotient rings $R = P/I$ for ideals $I\subset P$; 
    \item \verb|MPolyLocRing|: localizations $U^{-1}P$ of polynomial rings $P$ at multiplicative sets
      $U \subset P$. For the latter we support\footnote{We deliberately restrict ourselves to those types 
      actually needed for this particular publication.}
      \begin{itemize}
        \item \verb|MPolyPowersOfElement|: powers of a collection of elements $f_1,\dots, f_k \in P$;
        \item \verb|MPolyComplementOfPrimeIdeal|: complement of a prime ideal $Q \subset P$.
      \end{itemize}
    \item \verb|MPolyQuoLocRing|: localizations $U^{-1}R$ of quotients $R = P/I$, considered as a finite $P$-algebra,
      at multiplicative sets $U \subset P$ in polynomial rings as above.
  \end{itemize}
  Polynomial rings and their quotients have already been implemented in \Oscar 
  and its subsystems. Based on this work, the second named author developed the 
  localizations. 

  Note that we work exclusively with multiplicative sets in polynomial rings $P$. In particular, localizations 
  of quotient rings $R = P/I$ are realized as a localization of $R$, where we see $R$ as a finite $P$-module with a
  multiplication map. Successive localizations are automatically flattened to a single localization
  at products of multiplicative sets. There are various technical reasons 
  supporting this decision; one of them being that quotient rings are in general not factorial while 
  fast factorization of elements turns out to be key for performance of this implementation. 
  From a purely mathematical viewpoint, these choices for the 
  construction of localizations do not impose any relevant restrictions.

  \medskip
  A common feature shared by all these four types of rings is that they posses 
  \textit{coordinate functions}, 
  i.e. elements $x_1,\dots, x_n \in R$ such that every element $f \in R$ is a fraction 
  $f = \frac{p(x)}{q(x)}$ of polynomials $p, q \in \bbk[X_1,\dots, X_n]$ evaluated at the $x_i$. 
  Such coordinate functions are important as their images uniquely determine any homomorphism 
  \[
    \varphi \colon R \to Q, \quad \varphi(x_i) = g_i \in Q
  \]
  to an arbitrary commutative $\bbk$-algebra $Q$.
  By algebraic prolongation, an element $f=p/q$ as above is then taken to $\varphi(f) = p(g) \cdot q(g)^{-1} \in Q$.
  Based on this fact, homomorphisms from and to either of these four rings are supported in 
  \Oscar, together with their common operations such as computing kernels or preimages of ideals. 

\subsection{Functionality for ideals and modules}

  For each one of the four types of rings, as in Section \ref{sec:algebras}, we support the common functionality for ideals such as sums,
  intersections, radicals, minimal associated primes, primary
  decomposition, ideal quotients, saturation, and ideal 
  membership, and similar functionality for finitely generated free modules and subquotients. 

  The generic \Oscar framework for finitely generated (graded) modules represented as subquotients, and its instance for modules over polynomial rings, is due to J. B\"ohm and his student A. Dinges, see also \cite{Dinges}. The computational backend is \verb|Singular| \cite{Sin}.
  For the localizations and quotient rings 
  the second named author has partially picked up on the ideas communicated 
  in \cite{Posur18} for so-called \textit{computable rings}. Recall that a ring $R$ is computable 
  if there are algorithms which solve the following three problems: 
  \begin{enumerate}
    \item Given a morphism of free $R$-modules $\varphi \colon R^m \to R^n$, compute 
      a generating set $v_1,\dots, v_r \in R^m$ for the kernel of $\varphi$.
    \item Given a finitely generated submodule $I = \langle v_1,\dots, v_r\rangle \subset R^m$ 
      and an element $w \in R^m$, 
      decide whether $w \in I$. 
    \item Given an element $w \in I\subset R^m$ as above, find a set of elements $\lambda_1,\dots, \lambda_r \in R$ 
      such that $w = \sum_{i=1}^r \lambda_i \cdot v_i$.
  \end{enumerate}
  Now most algorithms involving ideals and modules in Commutative Algebra can be reduced to these 
  three building blocks. There is a layer of generic code for finitely presented 
  modules in \Oscar doing that.

  The key result in \cite{Posur18} is that for a computable ring 
  $R$, the localization at a multiplicative set $U \subset R$ is again computable if there is an algorithm for 
  the following problem:
  \begin{enumerate}
    \item[(4)] Given a finitely generated ideal $I \subset R$, decide whether $U \cap I$ is non-empty and 
      if so, produce a single element $x \in U \cap I$.
  \end{enumerate}
  Building on this result, the second named author 
  has extended the generic code base for finitely presented modules and ideals 
  so that it is now possible to have full functionality
  over localizations in \Oscar by only implementing one single method for a given new type of 
  multiplicative sets. 

  However, there are several problems for which a direct reduction to
  the above three methods does not yield a reasonable performance, or, even worse, 
  there is no known reduction to the above three methods. For
  these cases the \verb|julia| programming language allows us to clip
  in custom implementations based on the multiple dispatch system. 
  One example for a functionality which is not easily reduced to (1)-(3) above is primary 
  decomposition which is a backbone for
  modeling of algebraic cycles and intersection theory. 
  Another example where our actual implementations do not rely on (4) is localization 
  at powers of elements $U = \{f^k : k \in \NN_0\}$, $f \in R$. For these, 
  working with saturations of ideals and modules has proven to be more effective.

\subsection{Affine schemes}

  In general, an affine scheme in \Oscar is merely a wrapper for one of the above 
  types of rings 
  and a morphism of affine schemes is an arrow-reversing wrapper of the ring homomorphisms 
  for the pullback of regular functions. The abstract type for affine schemes and their 
  morphisms for which their minimal interfaces are defined, are \verb|AbsAffineScheme| and 
  \verb|AbsAffineSchemeMor|, respectively. 
  Minimal concrete implementations of these interfaces are realized for the 
  concrete types \verb|AffineScheme| and \verb|AffineSchemeMor|. An ordinary affine scheme 
  $X$ is constructed calling the function \verb|spec(R)| on a ring $R$ of the 
  above four types. The constructor for a morphism $X \to Y$ is 
  \verb|morphism(X, Y, [g_1, ..., g_n])| where $g_1,\dots, g_n \in R$ are the 
  pullbacks of the coordinate functions for $Y$ on $X = \Spec R$. 

  \subsubsection{Comparison of affine schemes}
  As already explained above, all of the four types of rings $R$ supported for 
  affine algebraic schemes have a governing polynomial ring $P = \bbk[X_1,\dots, X_n]$
  in the background, the \verb|ambient_coordinate_ring|. 
  Note that these arise by the nature of their \textit{implementation};
  mathematically these polynomial rings depend on a choice of the coordinate functions. 
  Even though this deviates from the strict mathematical definitions, we found it 
  reasonable and practical to allow for comparison of either two 
  affine schemes $X = \Spec R$ and $Y = \Spec R'$ based on their 
  embeddings $\Spec R \hookrightarrow \Spec P$ and $\Spec R' \hookrightarrow \Spec P$
  in their ``ambient spaces'' 
  whenever both $R$ and $R'$ are constructed using the same polynomial ring 
  $P$ as an object in memory.

  \subsubsection{Inheritance within the schemes framework}
  Even though the \verb|julia| programming language has no \textit{natural} 
  patterns for inheritance for its types, there are various possibilities 
  to mimic type-inheritance. For our schemes framework we decided to implement 
  the following simple pattern.
  \begin{itemize}
    \item Introduce an \textit{abstract type} \verb|MyAbstractType| 
      and declare an interface for that; 
      including all getters \verb|my_getter(::MyAbstractType)|, but 
      probably also more advanced functionality.
    \item Implement the generic getters via 
      \begin{verbatim}  my_getter(obj::MyAbstractType) = my_getter(underlying_concrete_instance(obj))\end{verbatim}
      and similarly for the other generic methods of the interface's functions.
    \item Introduce at least one \textit{minimal concrete type} 
      \verb|MyMinimalConcreteType <: MyAbstractType|
      together with a full implementation of the interface by methods for this 
      concrete type; i.e. \begin{verbatim}  my_getter(::MyMinimalConcreteType).\end{verbatim}
    \item For further, more specialized concrete types 
      \verb|NewConcreteType <: MyAbstractType|, either implement the interface 
      for \verb|MyAbstractType| directly by overwriting the respective methods, 
      or indirectly, by storing an instance of \verb|MyMinimalConcreteType| 
      inside \verb|NewConcreteType| and forwarding of its functionality via 
      implementation of 
      \verb|underlying_concrete_instance(::NewConcreteType)|. 
  \end{itemize}
  For instance, this allows us to introduce specialized types for 
  closed embeddings as \begin{verbatim}  ClosedEmbedding <: AbsAffineSchemeMor.\end{verbatim}
  Besides the usual data required for a morphism of affine schemes, 
  this also holds the information on the ideal defining the image in its ambient scheme. 
  While the multiple dispatch system of the \verb|julia| language allows for more 
  fine-grained handling of the functionality for such \verb|ClosedEmbedding|s, 
  forwarding a concrete instance of \verb|AffineSchemeMor| in the internals 
  enables the use of generic code written for \verb|AbsAffineSchemeMor| 
  and thus avoids code duplication. 

  Having described this pattern, and its importance for our implementation, we will 
  in the following only briefly comment on the abstract types for which an interface
  is defined, and the concrete types for their minimal implementations. For further 
  details, we refer the reader to the source code. 

\subsection{Covered schemes} 
Recall that a scheme over a ring $\bbk$ is a locally ringed space $(X,\mathcal{O}_X)$ which admits an open cover $(U_i)_{i \in I}$ such that each $(U_i,\mathcal{O}_X|U_i)\cong \Spec R_i$ as locally ringed spaces for some finitely generated $\bbk$-algebra $R_i$. 

In \Oscar, we start from finitely many affine patches $(U_i, \mathcal{O}_X|U_i)_{i\in I}=(\Spec R_i)_{i\in I}$ and isomorphisms $f_{i,j}$, called glue maps, between principal open 
subsets\footnote{We have partial support for gluings along more general open subsets,
but this is not required within the scope of this paper.}:
\begin{equation}\label{eqn:gluing}
  \Spec R_i = U_i \hookleftarrow U_{i, j} = 
  \xymatrix{
    \Spec R[h_j^{-1}] \ar@/^1pc/[r]^{f_{i, j}} &
    \Spec R_j[h_i^{-1}] \ar@/^1pc/[l]^{f_{j, i}}
  }
  = U_{j, i} \hookrightarrow U_j =\Spec R_j
\end{equation}
where we require for the glue maps that $f_{j,k} \circ f_{i,j} = f_{i,k}$ on the overlap $U_{i,j} \cap f_{i,j}^{-1}(U_{j,k})$.
On the mathematcal side, this means that we define
\begin{equation}
  \label{eqn:CoveredSchemes}
  Y := \bigsqcup_i U_i, \qquad X:=Y / \sim
\end{equation}
where for $x \in U_i$ and $y \in U_j$ we set $x \sim y$ if $f_{i,j}(x)=y$.
The structure sheaf $\mathcal{O}_X$ is the unique sheaf of rings such that the canonical projection $\pi \colon Y\to X$ is a morphism of locally ringed spaces.\\

Note that there is an implicit compromise in this construction: There is, for instance,
no \textit{single} version of $U_i \cap U_j$ for two affine charts $U_i$ and
$U_j$ of $X$. Instead, we are practically dealing with charts on $Y$.
In particular, $U_{i, j}$ and $U_{j,i}$ are not equal in memory (i.e. on $Y$)
but both \textit{represent} $U_i \cap U_j$.
However, introducing a new type of affine scheme
which would allow for different such representatives seemed excessive and would lead to inconsistencies with the design decisions made for affine schemes. Therefore
we treat $U_{i, j}$ and $U_{j, i}$ in the above sense as \textit{different}
schemes a priori, but we allow to ask for the identifying maps
$f_{ij}\colon U_{i, j} \to U_{j, i}$. This has further consequences down the road as we shall see in the next section on sheaves. 

In general, it is the user's responsibility to provide charts and gluings such that the resulting scheme is well
defined. For most practical applications, however, the construction of a scheme will
ensue from standard building blocks such as projective spaces or, 
more generally, toric varieties, and subschemes thereof, as well as
projective bundles, blow-ups, etc. In all such cases the gluings are provided automatically and 
separatedness of the scheme is guaranteed throughout the process.
\begin{remark}
The actual computation of automatically generated gluings can be rather expensive and,
given that the total number of gluings of a scheme $X$ is quadratic in the number $N$ of affine charts, 
should in general be avoided if possible. Therefore, most gluings are implemented lazy, 
i.e. they are computed only on request and then cached. 
\end{remark}

\medskip
The interface for gluings as in Equation \ref{eqn:gluing} is declared for the abstract type \verb|AbsGluing| and comprises
\verb|patches|, \verb|gluing_domains|, and \verb|gluing_morphisms| for the identifying maps in the middle.
The concrete type of any covered scheme in \Oscar is a subtype of the abstract type 
\verb|AbsCoveredScheme| and has at least one 
\verb|Covering|, its \verb|default_covering|. The latter then provides the affine \verb|patches| 
and their pairwise gluings. A basic implementation of the interface specified for \verb|AbsCoveredScheme|
is realized for the concrete type \verb|CoveredScheme|. Again, such basic implementations 
can be enriched with arbitrary additional data 
for specialized types and we have done so with our implementation of \verb|EllipticSurface|. Among 
other things, the latter also provides its \verb|weierstrass_model|, the associated elliptic 
curve over a function field via \verb|generic_fiber|, the \verb|trivial_lattice|, and the 
\verb|weierstrass_contraction| to project to the Weierstrass model.

\subsection{Coherent sheaves}

Recall that a \textit{presheaf} $\mathcal F$ on a topological space $X$ is a contravariant functor
which associates to every open subset $U \subset X$ an object $\mathcal F(U)$ and 
to every open inclusion of open sets $U \hookrightarrow V$ a \textit{restriction morphism}
\[
  \rho_{V, U} \colon \mathcal F(V) \to \mathcal F(U).
\]
A presheaf $\mathcal F$ is a \textit{sheaf}, if for any open $U \subset X$,
any covering
$\{U_i\}_{i \in I}$ of $U$, and any collection of sections
$s_i \in \mathcal F(U_i)$ with $s_i = s_j$ on $U_i \cap U_j$
there exists a unique global section $s \in \mathcal F(U)$ with $s_i = \rho_{U, U_i}(s)$.

A \textit{coherent sheaf} on a covered scheme $X$ as above is a sheaf of 
$\mathcal O_X$-modules which are locally finitely generated. On an affine scheme
$U = \Spec R$ any coherent sheaf $\mathcal F$ is uniquely determined by the finitely 
generated $R$-module $M = \mathcal F(U)$. 
\medskip

While the definition of these mathematical notions is quite straightforward, we are facing various constraints and challengens when attempting to model this in a computer:
\begin{itemize}
 \item We have to model a sheaf by a finite amount of data.
 \item We do not have direct access to the topology of $X$, only to $Y$ as in 
   (\ref{eqn:CoveredSchemes}) and the glue maps $f_{i,j}$.
\end{itemize}
In order to surmount these challenges, we model a coherent sheaf $\mathcal F$ on $X$ by its pullback along the canonical projection $Y \to X$. 
This forces us to make several compromises which lead to
a practically working substitute for the notion of a coherent sheaf, but at the expense that 
we slightly deviate from the strict mathematical concept. 

The first and probably most 
important compromise is that we actually only model \textit{presheaves}, 
i.e. we provide the functions $U\mapsto \mathcal F(U)$ on admissible 
open subsets $U$ and $U \supset V\mapsto \rho_{U,V}$ for admissible open inclusions. 
For the concrete implementations, we leave it as the programmer's responsibility 
to assure that they behave like honest sheaves.

Second, we restrict the form of open subsets $U \subset X$ which are recognized as \textit{admissible} for the production of the module $\mathcal F(U)$ for a coherent sheaf
$\mathcal F$; this will be detailed below.

The third compromise is inherent in our construction: the implicit identification of $U_i\subseteq Y$ with $\pi(U_i) \subseteq X$.
In consequence, objects which are
mathematically the same on $X$, appear from different representations on $Y$ by concrete data in memory:
For $U_i, U_j$ open in $X$ the module $\mathcal F(U_i \cap U_j)$ can be represented as
both, the localization $\mathcal F(U_{i,j})$ of the module $\mathcal F(U_i)$ and as localization $\mathcal{F}(U_{j,i})$ of $\mathcal F(U_j)$.
 Both representations are reconciled by giving an isomorphism $\rho_{ij}\colon \mathcal F(U_{i,j}) \to \mathcal F(U_{j,i})$ as part of the data.
%
In particular this holds for the structure sheaf $\OO_X$.

\medskip
The \textit{admissible} open subsets for a coherent sheaf $\mathcal F$ on a covered scheme $X$ are
\begin{itemize}
  \item the \verb|affine_charts| of $X$, i.e. the \verb|patches| $U_i$ of its \verb|default_covering| $(U_i)_{i \in I}$,
\end{itemize}
  and by recursive extension
\begin{itemize}
  \item any \verb|PrincipalOpenSubset| of an admissible open subset and
  \item any \verb|SimplifiedAffineScheme| of an admissible open subset.
\end{itemize}
The latter two types are special concrete subtypes of \verb|AbsAffineScheme|
which relate to an \verb|ambient_scheme|, resp. to an \verb|original_scheme|.
A \verb|PrincipalOpenSubset| arises as the complement of the zero locus
of a single element $h \in \OO_X(V)$, the \verb|complement_equation|,
on an affine scheme $V$.

\begin{remark}
The \verb|SimplifiedAffineScheme| was
introduced in order to eliminate superfluous variables whenever possible.
To this end, we apply Singular's \href{https://www.singular.uni-kl.de/Manual/4-3-0/sing_2139.htm}{elimpart} heuristic to the modulus $I$ of
a quotient ring $R = P/I$ for $V = \Spec R$ and wrap the result in an instance
of \verb|SimplifiedAffineScheme|, together with the identifying maps between both
representations of $\Spec R$.
Since the complexity of Buchberger's algorithm is doubly exponential
in the number of variables, this is a crucial step in keeping things computable;
especially in setups like iterated blowup constructions.
\end{remark}

With this recursively defined structure of a ``parent node'' given by
the ambient, resp. original scheme, the admissible open subsets of a covered
scheme $X$ form a forest (i.e. a union of trees in the graph-theoretic sense), whose
roots, the affine charts of $X$, are pairwise related by gluings.
This structure can then be used to indicate how to actually obtain
the module $\mathcal F(U)$ for admissible open subsets $U \subset X$ different from the
patches of the default covering: If $U$ is a \verb|PrincipalOpenSubset|, 
then $\mathcal F(V)$ is constructed as a localization of $\mathcal F(V)$ where 
$V$ is the \verb|ambient_scheme| of $U$. Similarly, in case $U$ is a \verb|SimplifiedAffineScheme|, 
the module $\mathcal F(U)$ is constructed via change of base rings from 
$\mathcal F(V)$ along the identifying morphisms with the \verb|original_scheme|
$V$. This pattern eventually recurses up until some patch is reached, on which 
a module for $\mathcal F$ is already known; be it by the original construction of 
$\mathcal F$, or a previous call.

Similarly, the tree structure is used to construct the restriction maps 
$\mathcal F(V) \to \mathcal F(U)$ for $U \subset V$ in $X$. If $U$ and 
$V$ come to lie in different affine charts, this involves the gluings 
on the root nodes of the respective two trees. 


\begin{remark}
The complement
$V$ of $\{xy - z^2 = 0\} \subset \PP^2$ is affine, but not contained in any of the standard affine charts of
$\PP^2$, seen as a covered scheme. To cover such cases, further types of admissible subsets may be added to \Oscar in the future.
\end{remark}

\medskip
The abstract type for which the user-facing interface for coherent sheaves is declared, 
is \verb|AbsCoherentSheaf|. 
The overall design of coherent sheaves in \Oscar is (and has to be) lazy and 
cached: No instance $\mathcal F(U)$ is computed, unless required for the definition 
of a coherent sheaf or requested by an explicit call. In general, the programmer is 
free to choose how the modules $\mathcal F(U)$ and their restriction morphisms 
are produced for their concrete subtypes of \verb|AbsCoherentSheaf|. 
While the constructor of the basic implementation requires 
given modules and identification maps 
over the affine charts of $X$ and their gluings, there are also special types like 
e.g. 
\verb|PullbackIdealSheaf|, the constructor of which 
requires only an \verb|AbsCoveredSchemeMorphism|
and an \verb|AbsIdealSheaf| on its codomain. Yet another type of ideal sheaf is 
\verb|PrimeIdealSheafFromChart| whose constructor only needs one prime ideal 
$I \subset \OO_X(U)$ on one non-empty admissible affine open $U \subset X$ 
as input and then extends to other charts from there. 
All of these types come with different patterns 
for their computation on admissible open subsets, realized by different sets of 
methods of the internal function \verb|produce_object|. See the source code for 
details and more, concrete examples. 

\subsection{Morphisms of covered schemes}

Let $X = \bigcup_{i \in I} U_i$ and $Y = \bigcup_{j\in J} V_j$ be covered schemes. In general,
a morphism $f \colon X \to Y$
is given by means of a \textit{refinement} $X = \bigcup_{i \in I, j\in J_i} U_{i, j}$ of the covering of $X$
and morphisms of affine schemes
\begin{equation}\label{eqn:morphism}
  f_{i, j} \colon U_{i, j} \to V_j,
\end{equation}
which glue naturally on the respective overlaps.
In \Oscar this information is held in a \verb|CoveringMorphism| which is in turn wrapped in some
concrete instance of \verb|AbsCoveredSchemeMorphism|.

For many setups this already provides convenient data structures to work with. For instance for
the blow-up $\pi \colon \mathrm{Bl}_{\mathcal I} X \to X$ of an ideal sheaf $\mathcal I$ on a
covered scheme $X$ the standard constructions of $\mathrm{Bl}_{\mathcal I}$ as a covered scheme
from $X = \bigcup_{i \in I} U_i$ leads to a covering
$\mathrm{Bl}_{\mathcal I} = \bigcup_{j \in J} V_j$ of $Y$ for which
one does not need a refinement to write down the local morphisms $V_j \to U_{i(j)}$ of affine
schemes for the projection map.

This example occurs in our implementation of elliptic surfaces 
when we construct a relatively minimal model $X$ as the minimal resolution of its (minimal) Weierstrass model $W$. The resolution $X \to W$ is a sequence of blowups and, in particular, an \verb|AbsCoveredSchemeMorphism|.

\subsection{Rational functions}
A \textit{variety} is a geometrically integral, separated scheme of finite type over a field $\bbk$. The corresponding type in \Oscar is \verb|AbsCoveredVariety|.
To a variety $X$ we can associate its field
of rational functions
\[
  k(X) = \mathrm{Quot}(\OO_X(U))
\]
which is isomorphic to the fraction field of the ring $\OO_X(U)$ of regular functions of any nonempty affine open subset $U \subset X$.
Consequently, an element $f \in k(X)$ comes with
many representatives: At least one for every affine chart of $X$.

The type for elements of $k(X)$ is \verb|VarietyFunctionFieldElem|.
Any such element can be asked to produce a representative on a given affine open subset
$V \subset X$ which is \textit{admissible} in the sense explained in the previous section
on coherent sheaves. Representatives are provided
as fractions of elements in the \verb|ambient_coordinate_ring|; a decision owed to the fact that
we do not support fraction fields of integral quotients of polynomial rings in \Oscar.

The implementation of rational functions on $X$ is lazy in the sense that representatives
are computed only in the moment where the user asks for one. At first, only one representative
$(U, f)$ with $f = \frac{p}{q}$ for one affine chart $U \subset X$ is stored in the data structure for $f$.
Once the user asks for a representative on an admissible open subset $V \subset X$,
a heuristic assessement of the presumed complexity of the production from already existing
representatives $(U', f')$ is made.
Such production involves evaluation of the gluing of $U'$ with $V$
which must be assumed to be computationally expensive in general.

\subsection{Birational Maps}
For morphisms like the sought for automorphism of an elliptically fibered K3 surface $f_{OY} \colon X_{OY} \to X_{OY}$ it is prohibitively expensive to compute a covering of $X_{OY}$ which is sufficiently fine-grained so that $f_{OY}$ can be described by local maps on its patches. However, this is not needed for our purposes such as computing the pushforward of divisors. It is therefore better to view $f_{OY}$ as a birational map. Let us recall the definition:\\

Let $X,Y$ be varieties, $U,U' \subseteq X$ and $V,V' \subseteq Y$
be non-empty open subsets. Two morphisms $\varphi\colon U \to V$ and
$\varphi'\colon U'\to V'$ are said to be \textit{equivalent} if they
agree on $U\cap U'$ (which is non-empty because $X$ is a variety). The
corresponding equivalence class is called a \textit{birational map} and it is
called \textit{dominant} if the image of any representative is dense in $Y$.\\

Within the scope of this paper, we work with dominant birational maps
$f \colon X \dashrightarrow Y$ 
which extend (uniquely) to a morphism $X \to Y$, that is, which have a trivial indeterminacy locus.
We nevertheless use birational morphisms and their representatives in the internals
of the implementation, as we shall see below.

A dominant birational map $f\colon X \dashrightarrow Y$ is uniquely determined 
by its pullback of function fields $k(Y) \to k(X)$ and this is also the
input data provided from our computations: Given the Weierstrass chart
$W_Y$ on an elliptic K3 surface $Y$ with coordinates $(t, x, y)$, the
computations from the $2$-neighbor step provide us with a Weierstrass
chart $W_X$ of $X$ and the pullbacks $ f^*(t),  f^*(x)$, and $ f^*(y)$
as rational functions on $W_X$. Since $K_X$ is nef and $X$ a (smooth)
surface, we know a priori that $f\colon X \to Y$ is an isomorphism
(although it is represented as a rational map).

\sloppy
All this is incorporated in a specialized
type for \verb|AbsCoveredSchemeMorphism| called \verb|MorphismFromRationalFunctions|. 
An approach similar to ours, but with a view towards applications in massively 
parallel setups, has also been pursued and implemented by 
B. Mirgain in \verb|Singular|; see \cite{Mirgain24}. 
Any
concrete instance of this type does
not attempt to compute a representation $\{f_{i, j}\colon U_{ij} \to V_j\}$ of itself as in Equation \ref{eqn:morphism}, unless requested by the user.
The functionality we need for such an automorphism --basically pullback and pushforward of divisors--
is mostly overwritten by methods which do not require the computation of a full
\verb|CoveringMorphism|. This is
crucial to the success of our computations, as realizing $f_{OY}$ via
a full-fledged \verb|CoveringMorphism| will probably not terminate in practice. 

\subsection{Specialization}
Recall that the base change of $X/R$ along a morphism $R \to k$, is defined by the fiber product
\[\begin{tikzcd}
X\times \Spec k \arrow[r,"\mathrm{red}"] \arrow[d] & X \arrow[d] \\
\Spec(k) \arrow[r]                   & \Spec(R).
\end{tikzcd}\]
\Oscar can compute the fiber product and the corresponding morphisms of 
\verb|AbsCoveredScheme|s by looking at the individual charts and their gluings.
By functoriality of the fiber product, local constructions patch together 
to give well-defined global objects.

In our application, we are facing the following caveat. The instances of 
\verb|AbsCoveredScheme|s that we are working with are \verb|EllipticSurface|s 
which hold significantly more information than just a covered scheme, and in positive
characteristic we would like to have the counterparts of this information available 
again. 
In particular, recall from Section \ref{sec:specialization} that we need access to the specialization map $\mathrm{sp}\colon \NS(X) \to \NS(X_k)$, 
so a ``forgetful'' reduction map which only has a \verb|CoveredScheme| in its codomain 
will not suffice. 

\sloppy
Now we could attempt to overwrite the method for reduction to positive
characteristic for an \verb|EllipticSurface|, but that would be technically 
very challenging, as we would have to dynamically adjust all cached objects such as, 
for instance, the Weierstrass models and projections, the basis of the \verb|algebraic_lattice|, etc.
Instead, we pursue a different approach. Given an \verb|EllipticSurface| $X$ over a ring $R$ of 
characteristic zero and a reduction map $R \to k$, we first specialize the
Weierstrass chart $W$ of $X$ to obtain the Weierstrass chart $W_k$ of
$X_k$ and construct the latter as an a-priori independent \verb|EllipticSurface|.
Via the identification
\[X \times \Spec k \hookleftarrow W_k \hookrightarrow X_k\]
we obtain a birational map
\[X \times \Spec k \dashrightarrow X_k\]
(which is in fact an isomorphism)
defined as a \verb|MorphismFromRationalFunctions|. Then we obtain the specialization map as the pushforward 
of the composition
\[X \overset{\mathrm{red}}{\longleftarrow} X \times \Spec k \overset{\cong}{\longrightarrow} X_k\]
On the true \verb|EllipticSurface|
$X_k$ we can now identify any given divisor as an element in $\NS(X_k)$ again 
and we get an identification of the generators for the algebraic lattices on both $X$ and $X_k$.
Similarly, we obtain $(f_{OY})_k$ as a \verb|MorphismFromRationalFunctions|  by specialization of its representation as a rational map $W \dashrightarrow W$ of the Weierstrass chart only.

Another advantage of this approach over the direct one is that data for $X_k$, such as the gluings between charts of $X_k$, can be computed directly in positive characteristic
rather than in characteristic $0$ (which is slower) and performing a base change.

\subsection{Algebraic cycles and divisors}

Following \cite{Fulton98}, an \textit{algebraic cycle} of dimension $k$ 
on a scheme $X$ is a formal
linear combination\footnote{Unless otherwise specified, we will assume integer coefficients throughout.}
\[
  C = \sum_{i=1}^N c_i \cdot W_i
\]
of irreducible subvarieties $W_i \subset X$ of dimension $k$. 
While in mathematical theory it comes at no expense to make such irreducibility 
assumptions on the summands of an algebraic cycle, this is not practical in an 
implementation; mostly due to the computational cost of primary decomposition.
Moreover, while it is theoretically possible to introduce a new concrete type 
of \verb|AbsCoveredScheme| for subschemes or -varieties of some ambient scheme, 
we have not taken this step, yet. Instead, we decide to model an algebraic cycle 
as a formal linear combination of ideal sheaves
\[
  C = \sum_{i=1}^N c_i \cdot \mathcal I_i.
\]
with $\mathcal I_i$ representing a subscheme $Y_i$ of $X$ which is not necessarily 
reduced or irreducible. If needed, the user can transform any such cycle 
following \cite[Section 1.5]{Fulton98} via \verb|irreducible_decomposition| 
into a linear combination
\[C = \sum_{\mathcal P \textnormal{ prime}} \left(\sum_i 
  c_i \cdot \operatorname{length}_{\OO_{X, \mathcal P}} (\OO_{X, \mathcal P}/(\mathcal I_i)_\mathcal P)\right) \cdot \mathcal{P}
\]
where 
the sum ranges over all \textit{prime} ideal sheaves $\mathcal{P}$, i.e. ideal sheaves
whose scheme theoretic zero loci are irreducible subvarieties of $X$ and 
$\OO_{X, \mathcal P}$ is the local ring at $\mathcal P$ 
(isomorphic to $(\OO_X(U))_{\mathcal P(U)}$ for any affine open $U \subset X$ 
is any affine chart with $\mathcal P(U) \neq \langle 1\rangle$).

The user interface for algebraic cycles on covered schemes in Oscar is
declared for the abstract type \verb|AbsAlgebraicCycle|. The basic implementation
is based on a dictionary, the keys of which are the \verb|components| $\mathcal I_i$, 
pointing to their coefficients $c_i$. 
If $X$ itself is equidimensional, a \textit{Weil divisor}
is an algebraic cycle of codimension one. They basically share their functionality with 
the algebraic cycles but they come with their own abstract type \verb|AbsWeilDivisor|. 

\medskip

For \textit{Cartier divisors} we slightly deviate from the common description in the literature 
and implement a \verb|CartierDivisor| as a formal sum of \verb|EffectiveCartierDivisor|s;
see e.g. the \href{https://stacks.math.columbia.edu/tag/01WQ}{Stacks project}. 
An \verb|EffectiveCartierDivisor| is internally based on an ideal sheaf $\mathcal I$ which is \textit{locally principal}
for some non-zerodivisor. 
In particular, any \verb|EffectiveCartierDivisor| can be asked for the \verb|trivializing_covering| and for every \verb|patch| $V$ of that covering the ideal 
$\mathcal I(V)$ is then of the desired form.

Within the scope of this paper Cartier divisors are naturally encountered as exceptional divisors 
from the resolution of singularities of the Weierstrass model of an elliptically fibered K3 surface. 
Every Cartier divisor can be transformed into a Weil divisor, but in general not the other way around.

\subsection{Pushforward and pullback of divisors}

Weil divisors or, more generally, algebraic cycles are covariant objects for \textit{proper} maps 
while Cartier divisors correspond -- up to rational equivalence -- to algebraic line bundles and 
are therefore of contravariant nature.
In \Oscar we have partial support for the associated \verb|pushforward| and \verb|pullback| functionality
for \verb|AbsCoveredSchemeMorphism|s. The
general pushforward of algebraic cycles for proper maps requires the computation of the degree of a field 
extension of fields of rational functions on varieties (cf. \cite[Section 1.4]{Fulton98}), a functionality 
which is not yet available in \Oscar.

In this paper we are mostly dealing with \textit{isomorphisms} of \textit{varieties} which are modeled as instances of
\verb|MorphismFromRationalFunctions|. In this setting pushforward and even pullback of algebraic cycles
is straightforward and can be reduced to mapping (non-closed) points of $X$.
\medskip

The implementation of 
the pushforward of an algebraic cycle $D$ along a \verb|MorphismFromRationalFunctions| $f\colon X \to Y$ is based 
on the following observation. Let $\mathcal P$ be a prime ideal sheaf on $X$. 
If $U \subset X$ is any affine chart of $X$ with $\mathcal P(U) \neq \langle 1 \rangle$, then the open subsets 
$U \cap f^{-1}(V_j)$ cover $U$ as $V_j$ ranges through the affine charts of $Y$. Hence, there
must be at least one $V = V_j$ for which $\mathcal P(U \cap f^{-1}(V)) \neq \langle 1\rangle$. 
In case $U \cap f^{-1}(V)$ is not principally open in $U$, we may restrict to some principally 
open subset $W \subset U \cap f^{-1}(V)$ with $\mathcal P(W) \neq \langle 1\rangle$.  
Once we found such a pair $(U, V)$, we can realize a local representative $f \colon U \cap f^{-1}(V) \to V$ 
of $f$, compute the scheme-theoretic image of the point $\mathcal P$ along this morphism of affine schemes 
and return the corresponding ideal sheaf on $Y$.

Given that the rational map $f$ in our application is huge, the practical difficulty is to find a suitable pair $(U,V)$ such that the computation of the image $f(\mathcal P(U))$ terminates.
Thus, we proceed as follows:
\begin{enumerate}
  \item Assume that the divisor is of the form $D = \sum_i d_i \cdot \mathcal P_i$ 
    for \textit{prime} ideal sheaves $\mathcal P_i$
    and reduce by linearity to the pushforward of a prime component
    $\mathcal P = \mathcal P_i$.
  \item Find one affine chart $U \subseteq X$ in the \verb|default_covering| of $X$ where $\mathcal P$ is visible, i.e.
    where $\mathcal P(U) \neq \langle 1\rangle$.
  \item Iterate through the \verb|patches| $V_j\subseteq Y$ of the \verb|default_covering| of the codomain $Y$ of $f$
    and compute the pullbacks of the coordinate functions $y_1,\dots, y_k$ for the ring $\OO_Y(V_j)$ 
    as rational functions on $U$ for all $V_j$. 
  \item Sort the codomain patches $V_j$ based on the complexity of the rational functions for $f$ on the pairs $(U, V_j)$ 
    and iterate through the $V_j$ in that order.
    Realize $f$ on \textit{some} principal open subset $W = D(h) \subset U \cap f^{-1}(V)$ 
    as a morphism of affine schemes
    \[
      f \colon W \to V.
    \]
    The element $h \in \OO_X(U)$ arises from the denominators of the fractions 
    for the local representation of $f$ on this 
    pair of patches $(U, V)$.
    If $\mathcal P(W) \neq \langle 1\rangle$, compute the scheme theoretic image of this point and return
    the associated prime ideal sheaf on $Y$.
  \item If the previous step fails, iterate through the $V_j$ again. For every $V = V_j$
    compute the maximal Zariski open subset $W' \subset U$ which allows for a realization of 
    $f \colon W' \to V$ as a morphism of schemes.
    Such $W'$ can not be expected to be principally open in $U$, but it can be covered by principal 
    open subsets. Eventually we must hit some $W'$ for which $\mathcal P(W') \neq \langle 1 \rangle$ 
    and we can finish as above.
\end{enumerate}
Correctness of this algorithm then follows from the realizability of $f$ as a morphism of covered schemes. 

\medskip
Pullback of Weil divisors along an isomorphism $f \colon X \to Y$ given as a 
\verb|MorphismFromRationalFunctions| is implemented similarly. Depending on the use case, sometimes the 
pushforward is faster, sometimes the pullback. 

\medskip
For an automorphism $f$ of an \verb|EllipticSurface| $X$ we have implemented an even more specialized 
method for the pushforward of the known \verb|algebraic_lattice| of $X$. 
For this we make a pre-assessement 
to determine the nature of the image of any component.
This is achieved by 
composing $f$ with the projection $\pi \colon X \to \PP^1$ of the elliptic fibration and pushforward 
of the components to $\PP^1$: If their dimension is one, the respective component is mapped to a 
(multi-) section in $X$, if the dimension is zero, it is taken to a fiber component.
This information 
can be used to significantly narrow down the affine charts $V$ in the codomain of $f$ for which one 
has to create local realizations of $f$ and we get significant speedups compared to the generic method 
for pushforward.

\subsection{Computation of intersection numbers}
On an elliptically fibered K3 surface $X$ computations of intersection numbers are important for the Gram-matrix
of the \verb|algebraic_lattice| and for representing the numerical equivalence class of an
arbitrary divisor $D$ as a linear combination of a chosen generating set of the N\'eron-Severi group
of $X$. Let us first take a look in greater generality:\\

Following \cite[Definition 2.4.1 and Section 2.5]{Fulton98}, the natural intersection product on a variety $X$ 
over $\bbk$ is a is defined for 
a Cartier divisor\footnote{In fact, Fulton introduces a \textit{pseudo-divisor} to address various technical difficulties; we shall not go into these details here.} $C$ via the first Chern class of the associated line bundle,
seen as a map 
\[
  c_1(C) \colon Z_k(X) \to A_{k-1}(X), \quad D \mapsto c_1(C)\cap D
\]
from the $k$-cycles on $X$ to the \textit{Chow group} of $(k-1)$-cycles on $X$.
The \textit{integral} of a zero-dimensional cycle $D = \sum_{\mathcal P} c_\mathcal P \cdot \mathcal P \in Z_0(X)$ is 
\begin{equation}
  \label{eqn:DefinitionIntegral}
  \int_X D = \sum_{\mathcal P} c_\mathcal P \cdot [k(P):\bbk]
\end{equation}
where we assume all $\mathcal P$ to be prime. For an (equidimensional) surface $X$ 
which is proper over $\bbk$ one can then 
define the intersection number of a Cartier and a Weil divisor $C$ and $D$ as 
\[
  C.D = \int_X c_1(C) \cap D.
\]
Whenever $X$ is smooth, every Weil divisor is actually Cartier, so that the definition 
extends to pairs of Weil divisors. In this case, it is customary to declare the intersection 
of two Weil divisors to be this number, provided they are in sufficiently general position 
relative to each other; cf. \cite[Chapter V, Theorem 1.1]{Hartshorne77}.
Since intersection numbers do not change under linear equivalence, this condition can always be achieved by a moving lemma, in theory at least.\\

\medskip
While the theory is straightforward and well developed, one quickly reveals severe limitations of 
the above definitions when it comes to their practical feasibility in applications like ours. 
One of the main bottlenecks is again the requirement to have a Weil divisor as a linear combination 
of \textit{irreducible} components. Fortunately, the computation of such a representation of a given 
cycle of dimension zero can be avoided if one is only interested in the integral of that cycle. We will describe 
this in more detail below. 

Another instance where irreducibility of components is relevant in 
intersection theory is \textit{self-intersection} of Cartier divisors, or, more generally, 
dealing with pairs of a Cartier and a Weyl divisor $C$ and $D$ 
where some irreducible component of $D$ happens to 
be properly contained in the \textit{support} $|C|$ of the Cartier divisor $C$. Usually 
one employs some kind of ``moving lemma'' to resolve such situations. While we do have the 
infrastructure necessary to implement such moving lemmas, we shall not need them for the purpose 
of this article: On our elliptically fibered K3 surfaces all relevant divisors in the algebraic 
lattice are already given as irreducible curves and their self intersection is completely determined 
by their genus via the adjunction formula and the Riemann-Roch theorem for curves. 
Thus we deliberately hardcode such excess- and self-intersections 
whenever encountered in our specialized methods for an \verb|EllipticSurface|.

\medskip

In order to bring the definition of the integral, Equation (\ref{eqn:DefinitionIntegral}), in a 
workable form for computation of intersection numbers on surfaces, we start with the following 
observations. 
Fix a component $\mathcal P$ of $D$ and let $P = V(\mathcal P) \subset X$ be its zero 
locus; choose an affine chart $U \subset X$ for which $\mathcal P(U) \neq \langle 1\rangle$. 
Then the degree of the field extension 
\[
  [k(P):\bbk] = \dim_\bbk \OO_X(U)/\mathcal P(U)
\]
is nothing but a vector space dimension of a quotient ring. Now if $\mathcal Q$ is $\mathcal P$-primary, 
then $[\mathcal Q] = n_{\mathcal P} \cdot [\mathcal P]$ with 
\[
  n_{\mathcal P} = \mathrm{length}_{\OO_{X, P}} \OO_{X, P}/\mathcal Q
\]
the \textit{length} of the Artinian quotient ring $\OO_{X, P}/\mathcal Q$; 
cf. \cite[Chapter 1, Section 1.5]{Fulton98}. Since we are working with zero dimensional 
cycles, it is easy to see that then 
\[
  n_P \cdot [k(P): \bbk] = \dim_\bbk \OO_X(U)/\mathcal Q(U)
\]
so that we can compute the full contribution of $\mathcal Q$ to the integral directly without 
extracting the prime component $\mathcal P$. 

For an arbitrary zero-dimensional ideal sheaf $\mathcal I$ on $X$ we observe that according to 
the Chinese Remainder Theorem the quotient $\mathcal O_X(U)/\mathcal I(U)$ decomposes as a direct sum 
\[
  \mathcal O_X(U)/\mathcal I(U) \cong \bigoplus_i \OO_X(U)/Q_i
\]
for the primary components $Q_i$ of $\mathcal I(U)$ in every affine chart $U$ of $X$. Therefore 
the contribution to an integral $\int_X 1 \cdot \mathcal I$ from the components visible in a 
fixed affine chart $U \subset X$ is 
\[
  \dim_\bbk \OO_X(U)/\mathcal I(U).
\]
This boils down to the computation of one single Gr\"obner basis rather than primary decompositions 
and computations of lengths.

At this point the main technical challenge is to avoid overcounting: One can compute contributions to 
the integral $\int_X D$ from every one of the affine charts $U$ of $X$, but in general the irreducible 
components $\mathcal P$ of $D$ will be visible in more than just one chart. 

Without loss of generality we may assume that $D = 1\cdot \mathcal I$ is given by a single ideal sheaf
$\mathcal I$ on $X$. In principle, it
is possible to compute local primary decompositions of $\mathcal I(U)$
and proceed with matching of the components via the gluing maps of either
two charts. However, this turns out to be computationally rather involved
with inferior performance: Not only is primary decomposition expensive, but also the computation
of gluings for the matching of local components should be avoided. 

To circumvent this problem, we write $X$ as a finite disjoint union of locally closed subsets $C_i$
which are contained as closed affine subschemes in the affine charts $U_i$ of $X$. 
In other words, we view
\[X= \bigsqcup_i C_i\]
as a constructible set.
Note that there is no natural choice for such decompositions, but
for most constructions it is possible to write them down. A common example for this is
the decomposition
\[
  \PP^n
  = \IA^n \sqcup \PP^{n-1}
  = \IA^n \sqcup \IA^{n-1}
  \sqcup \cdots
  \sqcup \IA^1
  \sqcup \IA^0
\]
of projective space into the disjoint union of affine spaces, starting with the horizon
$\PP^{n-1} \subset \PP^n$ at infinity and its complement $\IA^n$ and the inductive continuation
of this procedure with $\PP^{n-1}$. Since the preimage of a constructible set under a morphism of schemes is constructible, we can pull back such a decomposition along sequences of blowups, forming total spaces of projective bundles, taking subschemes or finite covers and refine the decomposition further if necessary. This is done automatically in \Oscar.

Technically, this is modeled by what we colloquially call \verb|decomposition_info|. For a given \verb|Covering| $(U_i)_{i\in I}$ of a covered scheme $X$ the decomposition info is a dictionary which holds for every \verb|patch| $U_i$ of the covering a list of elements
$f_{1,i},\dots, f_{k_i,i} \in \OO_X(U)$. These elements are chosen so that their common zero loci $C_i:= V(f_{1,i},\dots, f_{k_i,i}) \subset U_i \subset X$ provide a decomposition of $X$ as a \textit{disjoint union} of
locally closed subsets.

For the computation of an integral the \verb|decomposition_info| comes in precisely to avoid overcounting: 
The contribution 
to $\int_X 1\cdot \mathcal I$ from a single chart $U_i$ is computed as follows.
We
throw away all components of $I = \mathcal I(U)$ whose support does not lie in $C_i$. This only needs saturation and
can be done without computing a primary decomposition.
We obtain an ideal $J_i \subset \OO_X(U_i)$ for the remaining components and the
contribution to the integral from this chart $U_i$
then amounts to $\dim_\bbk \OO_X(U_i)/J_i$.

\section{Appendix}
\begin{landscape}
\begin{figure}
  \label{fig:LiftingsToCharZero}
\begin{equation*}
  \begin{matrix}
    \hline
    \left(
    \begin{matrix}
      (\frac{1}{2}\cdot \alpha^7 + \frac{5}{2}\cdot \alpha^3 - \alpha^2 - \alpha)\cdot t^4 \\
      + (\frac{1}{2}\cdot \alpha^7 + \frac{1}{2}\cdot \alpha^5 + \alpha^4 + \frac{7}{2}\cdot \alpha^3 - \alpha^2 - \frac{1}{2}\cdot \alpha)\cdot t^3 \\
      + (-\frac{1}{2}\cdot \alpha^6 + \alpha^4 + 2\cdot \alpha^3 - \frac{1}{2}\cdot \alpha^2)\cdot t^2 \\
      + (\frac{1}{2}\cdot \alpha^7 - \frac{1}{2}\cdot \alpha^6 - \frac{1}{2}\cdot \alpha^5 - \frac{1}{2}\cdot \alpha^4 + \frac{7}{2}\cdot \alpha^3 - \frac{1}{2}\cdot \alpha^2 + \frac{1}{2}\cdot \alpha - \frac{1}{2})\cdot t \\
      + \frac{1}{2}\cdot \alpha^7 - \frac{1}{2}\cdot \alpha^4 + \frac{5}{2}\cdot \alpha^3 + \alpha - \frac{1}{2}
    \end{matrix}\right., &
    \left.
    \begin{matrix}
      (\frac{1}{2}\cdot \alpha^6 + \frac{1}{2}\cdot \alpha^5 + \frac{1}{2}\cdot \alpha^4 + \frac{3}{2}\cdot \alpha^2 - \frac{1}{2}\cdot \alpha - \frac{1}{2})\cdot t^6 \\
      + (-\alpha^7 + \alpha^5 + 3\cdot \alpha^4 - \alpha^3 + 3\cdot \alpha^2 + \alpha)\cdot t^5\\
      + (-\frac{3}{2}\cdot \alpha^7 - 2\cdot \alpha^6 - 3\cdot \alpha^5 - \frac{3}{2}\cdot \alpha^3 + 4\cdot \alpha^2)\cdot t^4 \\
      + (\alpha^7 - \frac{1}{2}\cdot \alpha^6 - 5\cdot \alpha^5 - \frac{15}{2}\cdot \alpha^4 + \alpha^3 + \frac{9}{2}\cdot \alpha^2 - \alpha - \frac{1}{2})\cdot t^3\\
      + (\frac{3}{2}\cdot \alpha^7 + \frac{5}{2}\cdot \alpha^6 - \frac{9}{2}\cdot \alpha^4 - \frac{3}{2}\cdot \alpha^3 + \frac{7}{2}\cdot \alpha^2 + \frac{1}{2})\cdot t^2\\
      + (\frac{3}{2}\cdot \alpha^6 + 2\cdot \alpha^5 + \frac{3}{2}\cdot \alpha^4 - 2\cdot \alpha^3 + \frac{3}{2}\cdot \alpha^2 + \frac{3}{2})\cdot t \\
      + \frac{1}{2}\cdot \alpha^5 + \alpha^4 + \alpha^3 + \frac{1}{2}\cdot \alpha + 1
    \end{matrix}\right)
    \\
    \hline

    \left(
    \begin{matrix}
      -t^4 + (\frac{1}{2}\cdot \alpha^6 + \frac{5}{2}\cdot \alpha^2 - 1)\cdot t^2 \\
      + \frac{1}{2}\cdot \alpha^6 + \frac{5}{2}\cdot \alpha^2
    \end{matrix}, \right. & 
    \left.
    \begin{matrix}
      (-\frac{1}{2}\cdot \alpha^7 - \frac{1}{2}\cdot \alpha^5 - \frac{5}{2}\cdot \alpha^3 - \frac{5}{2}\cdot \alpha)\cdot t^5 \\
      + (\frac{1}{2}\cdot \alpha^7 - \frac{1}{2}\cdot \alpha^5 + \frac{7}{2}\cdot \alpha^3 - \frac{7}{2}\cdot \alpha)\cdot t^3 
      + (\alpha^7 + 6\cdot \alpha^3 - \alpha)\cdot t
    \end{matrix} \right)
    \\
    \hline

    \left(
      t^4 - 1, 
      \right.
      &
      \left.
      \begin{matrix}
        (-\frac{1}{2}\cdot \alpha^6 - \frac{7}{2}\cdot \alpha^2)\cdot t^4 + \frac{1}{2}\cdot \alpha^6 + \frac{7}{2}\cdot \alpha^2
      \end{matrix} 
      \right)
      \\
    \hline

    \left(
    \begin{matrix}
      (\alpha^7 + \alpha^6 + \frac{1}{2}\cdot \alpha^5 + 6\cdot \alpha^3 + 6\cdot \alpha^2 + \frac{5}{2}\cdot \alpha)\cdot t^4 \\
      + (-\frac{7}{2}\cdot \alpha^7 - \alpha^6 + \frac{1}{2}\cdot \alpha^5 + \alpha^4 - \frac{41}{2}\cdot \alpha^3 - 6\cdot \alpha^2 + \frac{7}{2}\cdot \alpha + 6)\cdot t^3 \\
      + (-\frac{5}{2}\cdot \alpha^6 - 2\cdot \alpha^5 - \alpha^4 - \frac{29}{2}\cdot \alpha^2 - 12\cdot \alpha - 6)\cdot t^2 \\
      + (\frac{7}{2}\cdot \alpha^7 + \frac{5}{2}\cdot \alpha^6 + \frac{1}{2}\cdot \alpha^5 - \frac{1}{2}\cdot \alpha^4 + \frac{41}{2}\cdot \alpha^3 + \frac{29}{2}\cdot \alpha^2 + \frac{7}{2}\cdot \alpha - \frac{5}{2})\cdot t \\
      - \alpha^7 + \frac{1}{2}\cdot \alpha^5 + \frac{1}{2}\cdot \alpha^4 - 6\cdot \alpha^3 + \frac{5}{2}\cdot \alpha + \frac{5}{2}
    \end{matrix},\right. &
    \left.
    \begin{matrix}
      (\frac{7}{2}\cdot \alpha^7 + \frac{3}{2}\cdot \alpha^6 - \frac{1}{2}\cdot \alpha^4 + \frac{41}{2}\cdot \alpha^3 + \frac{17}{2}\cdot \alpha^2 - \frac{7}{2})\cdot t^6 \\
      + (-5\cdot \alpha^7 + 3\cdot \alpha^6 + 5\cdot \alpha^5 + 3\cdot \alpha^4 - 29\cdot \alpha^3 + 18\cdot \alpha^2 + 29\cdot \alpha + 18)\cdot t^5 \\
      + (-18\cdot \alpha^7 - 16\cdot \alpha^6 - \frac{15}{2}\cdot \alpha^5 - 105\cdot \alpha^3 - 94\cdot \alpha^2 - \frac{87}{2}\cdot \alpha)\cdot t^4 \\
      + (29\cdot \alpha^7 + \frac{15}{2}\cdot \alpha^6 - 5\cdot \alpha^5 - \frac{15}{2}\cdot \alpha^4 + 169\cdot \alpha^3 + \frac{89}{2}\cdot \alpha^2 - 29\cdot \alpha - \frac{89}{2})\cdot t^3 \\
      + (\frac{23}{2}\cdot \alpha^6 + \frac{21}{2}\cdot \alpha^5 + \frac{9}{2}\cdot \alpha^4 + \frac{133}{2}\cdot \alpha^2 + \frac{123}{2}\cdot \alpha + \frac{55}{2})\cdot t^2 \\
      + (-12\cdot \alpha^7 - \frac{15}{2}\cdot \alpha^6 - 2\cdot \alpha^5 + \frac{3}{2}\cdot \alpha^4 - 70\cdot \alpha^3 - \frac{87}{2}\cdot \alpha^2 - 12\cdot \alpha + \frac{15}{2})\cdot t \\
      + \frac{5}{2}\cdot \alpha^7 - \alpha^5 - \alpha^4 + \frac{29}{2}\cdot \alpha^3 - 6\cdot \alpha - 5
    \end{matrix} \right)
    \\
    \hline

    \left(
    \begin{matrix}
      (\frac{1}{2}\cdot \alpha^4 + \frac{1}{2})\cdot t^4 
      + (-\frac{1}{2}\cdot \alpha^4 + \alpha^2 - \frac{1}{2})\cdot t^2 
      - \alpha^2
    \end{matrix}, \right. &
    \left.
    \begin{matrix} 
      (-\frac{1}{2}\cdot \alpha^6 - \frac{3}{2}\cdot \alpha^2)\cdot t^6 
      + (\frac{1}{2}\cdot \alpha^6 - \alpha^4 + \frac{1}{2}\cdot \alpha^2)\cdot t^4 \\
      + (\frac{3}{2}\cdot \alpha^4 + \frac{1}{2})\cdot t^2 
      - \frac{1}{2}\cdot \alpha^4 + \alpha^2 - \frac{1}{2}
    \end{matrix}
    \right)
    \\
    \hline
    
    \left(
    \begin{matrix}
      (\frac{1}{2}\cdot \alpha^7 + \frac{5}{2}\cdot \alpha^3 + \alpha^2 - \alpha)\cdot t^4 \\
      + (\frac{1}{2}\cdot \alpha^7 - \frac{1}{2}\cdot \alpha^6 + \frac{1}{2}\cdot \alpha^5 + \frac{1}{2}\cdot \alpha^4 + \frac{3}{2}\cdot \alpha^3 - \frac{1}{2}\cdot \alpha^2 + \frac{3}{2}\cdot \alpha + \frac{1}{2})\cdot t^3 \\
      + (-\frac{1}{2}\cdot \alpha^6 + \alpha^4 - 2\cdot \alpha^3 - \frac{1}{2}\cdot \alpha^2)\cdot t^2 \\
      + (\frac{1}{2}\cdot \alpha^7 - \frac{1}{2}\cdot \alpha^5 + \alpha^4 + \frac{3}{2}\cdot \alpha^3 + \alpha^2 - \frac{3}{2}\cdot \alpha)\cdot t \\
      + \frac{1}{2}\cdot \alpha^7 + \frac{1}{2}\cdot \alpha^4 + \frac{5}{2}\cdot \alpha^3 + \alpha + \frac{1}{2}
    \end{matrix}, \right. &
    \left.
    \begin{matrix}
      (-\frac{1}{2}\cdot \alpha^7 + \frac{1}{2}\cdot \alpha^6 - \frac{1}{2}\cdot \alpha^4 - \frac{3}{2}\cdot \alpha^3 + \frac{3}{2}\cdot \alpha^2 + \frac{1}{2})\cdot t^6 \\
      + (-\alpha^7 + \alpha^5 - 3\cdot \alpha^4 - \alpha^3 - 3\cdot \alpha^2 + \alpha)\cdot t^5 \\
      + (-\frac{3}{2}\cdot \alpha^7 + 3\cdot \alpha^5 - 7\cdot \alpha^4 - \frac{3}{2}\cdot \alpha^3 - 3)\cdot t^4 \\
      + (-\alpha^7 - \frac{1}{2}\cdot \alpha^6 + 5\cdot \alpha^5 - \frac{15}{2}\cdot \alpha^4 - \alpha^3 + \frac{9}{2}\cdot \alpha^2 + \alpha - \frac{1}{2})\cdot t^3 \\
      + (-\frac{5}{2}\cdot \alpha^6 + \frac{9}{2}\cdot \alpha^5 - \frac{9}{2}\cdot \alpha^4 - \frac{7}{2}\cdot \alpha^2 + \frac{3}{2}\cdot \alpha + \frac{1}{2})\cdot t^2 \\
      + (-\frac{3}{2}\cdot \alpha^6 + 2\cdot \alpha^5 - \frac{3}{2}\cdot \alpha^4 - 2\cdot \alpha^3 - \frac{3}{2}\cdot \alpha^2 - \frac{3}{2})\cdot t \\
      + \frac{1}{2}\cdot \alpha^5 - \alpha^3 + 2\cdot \alpha^2 + \frac{1}{2}\cdot \alpha
    \end{matrix} 
    \right)
    \\
    \hline

    \left(
    \begin{matrix}
      (-\alpha^7 - \alpha^6 - \frac{1}{2}\cdot \alpha^5 - 6\cdot \alpha^3 - 6\cdot \alpha^2 - \frac{5}{2}\cdot \alpha)\cdot t^4 \\
      + (-\frac{7}{2}\cdot \alpha^7 - \frac{5}{2}\cdot \alpha^6 - \frac{1}{2}\cdot \alpha^5 + \frac{1}{2}\cdot \alpha^4 - \frac{41}{2}\cdot \alpha^3 - \frac{29}{2}\cdot \alpha^2 - \frac{7}{2}\cdot \alpha + \frac{5}{2})\cdot t^3 \\
      + (-5\cdot \alpha^7 - \frac{5}{2}\cdot \alpha^6 + \alpha^4 - 29\cdot \alpha^3 - \frac{29}{2}\cdot \alpha^2 + 6)\cdot t^2 \\
      + (-\frac{7}{2}\cdot \alpha^7 - \alpha^6 + \frac{1}{2}\cdot \alpha^5 + \alpha^4 - \frac{41}{2}\cdot \alpha^3 - 6\cdot \alpha^2 + \frac{7}{2}\cdot \alpha + 6)\cdot t \\
      - \alpha^7 + \frac{1}{2}\cdot \alpha^5 + \frac{1}{2}\cdot \alpha^4 - 6\cdot \alpha^3 + \frac{5}{2}\cdot \alpha + \frac{5}{2}
    \end{matrix}, \right.& 
    \left.
    \begin{matrix}
      (\frac{3}{2}\cdot \alpha^6 + \frac{3}{2}\cdot \alpha^5 + \frac{1}{2}\cdot \alpha^4 + \frac{17}{2}\cdot \alpha^2 + \frac{17}{2}\cdot \alpha + \frac{7}{2})\cdot t^6 \\
      + (5\cdot \alpha^7 + \frac{15}{2}\cdot \alpha^6 + 5\cdot \alpha^5 + \frac{3}{2}\cdot \alpha^4 + 29\cdot \alpha^3 + \frac{87}{2}\cdot \alpha^2 + 29\cdot \alpha + \frac{15}{2})\cdot t^5 \\
      + (18\cdot \alpha^7 + 16\cdot \alpha^6 + \frac{15}{2}\cdot \alpha^5 + 105\cdot \alpha^3 + 94\cdot \alpha^2 + \frac{87}{2}\cdot \alpha)\cdot t^4 \\
      + (29\cdot \alpha^7 + \frac{37}{2}\cdot \alpha^6 + 5\cdot \alpha^5 - \frac{7}{2}\cdot \alpha^4 + 169\cdot \alpha^3 + \frac{215}{2}\cdot \alpha^2 + 29\cdot \alpha - \frac{37}{2})\cdot t^3 \\
      + (\frac{51}{2}\cdot \alpha^7 + \frac{23}{2}\cdot \alpha^6 - \frac{9}{2}\cdot \alpha^4 + \frac{297}{2}\cdot \alpha^3 + \frac{133}{2}\cdot \alpha^2 - \frac{55}{2})\cdot t^2 \\
      + (12\cdot \alpha^7 + 3\cdot \alpha^6 - 2\cdot \alpha^5 - 3\cdot \alpha^4 + 70\cdot \alpha^3 + 18\cdot \alpha^2 - 12\cdot \alpha - 18)\cdot t \\
      + \frac{5}{2}\cdot \alpha^7 - \alpha^5 - \alpha^4 + \frac{29}{2}\cdot \alpha^3 - 6\cdot \alpha - 5
    \end{matrix}\right)
    \\
    \hline
    
    \left(
    \begin{matrix}
      (\frac{1}{2}\cdot \alpha^7 + \frac{5}{2}\cdot \alpha^3 - \alpha^2 - \alpha)\cdot t^4 \\
      + (\frac{1}{2}\cdot \alpha^7 - \frac{1}{2}\cdot \alpha^6 - \frac{1}{2}\cdot \alpha^5 - \frac{1}{2}\cdot \alpha^4 + \frac{7}{2}\cdot \alpha^3 - \frac{1}{2}\cdot \alpha^2 + \frac{1}{2}\cdot \alpha - \frac{1}{2})\cdot t^3 \\
      + (\frac{1}{2}\cdot \alpha^6 + \alpha^5 + \alpha^4 + \frac{1}{2}\cdot \alpha^2 + \alpha)\cdot t^2 \\
      + (-\frac{1}{2}\cdot \alpha^7 - \frac{1}{2}\cdot \alpha^5 - \alpha^4 - \frac{7}{2}\cdot \alpha^3 + \alpha^2 + \frac{1}{2}\cdot \alpha)\cdot t \\
      - \frac{1}{2}\cdot \alpha^7 + \frac{1}{2}\cdot \alpha^4 - \frac{5}{2}\cdot \alpha^3 - \alpha + \frac{1}{2}
    \end{matrix}, \right.&
    \left.
    \begin{matrix}
      (\frac{1}{2}\cdot \alpha^6 + \frac{1}{2}\cdot \alpha^5 + \frac{1}{2}\cdot \alpha^4 + \frac{3}{2}\cdot \alpha^2 - \frac{1}{2}\cdot \alpha - \frac{1}{2})\cdot t^6 \\
      + (-2\cdot \alpha^5 - 3\cdot \alpha^4 - 2\cdot \alpha^3 + 3\cdot \alpha^2)\cdot t^5 \\
      + (-\frac{3}{2}\cdot \alpha^7 + 3\cdot \alpha^5 + 7\cdot \alpha^4 - \frac{3}{2}\cdot \alpha^3 + 3)\cdot t^4 \\
      + (2\cdot \alpha^7 + \frac{1}{2}\cdot \alpha^6 - 2\cdot \alpha^5 - \frac{15}{2}\cdot \alpha^4 - \frac{9}{2}\cdot \alpha^2 - \frac{1}{2})\cdot t^3 \\
      + (-\frac{3}{2}\cdot \alpha^7 - \frac{5}{2}\cdot \alpha^6 + \frac{9}{2}\cdot \alpha^4 + \frac{3}{2}\cdot \alpha^3 - \frac{7}{2}\cdot \alpha^2 - \frac{1}{2})\cdot t^2 \\
      + (\alpha^7 + \frac{3}{2}\cdot \alpha^6 + \alpha^5 - \frac{3}{2}\cdot \alpha^4 + \alpha^3 + \frac{3}{2}\cdot \alpha^2 + \alpha - \frac{3}{2})\cdot t \\
      - \frac{1}{2}\cdot \alpha^5 + \alpha^3 + 2\cdot \alpha^2 - \frac{1}{2}\cdot \alpha 
    \end{matrix}
    \right)
  \end{matrix}
\end{equation*}
  \caption{Liftings of the rational points in (\ref{eqn:RationalPointsPosChar}) on $\overline E_1$ from Section \ref{sec:FindingSections} to characteristic zero}
\end{figure}
\end{landscape}

\printbibliography
\end{document}